\documentstyle[12pt,amssymb]{amsart}
\renewcommand{\marginpar}[1]{}
\catcode`\@=12

\def\Empty{}
\newcommand\oplabel[1]{
  \def\OpArg{#1} \ifx \OpArg\Empty {} \else
  	\label{#1}
  \fi}
		
\long\def\realfig#1#2#3#4{
\begin{figure}[htbp]
\centerline{\psfig{figure=#2,width=#4}}
\caption[#1]{#3}
\oplabel{#1}
\end{figure}}

\newcommand{\comm}[1]{}
    
\input{psfig}
\newtheorem{thm}{Theorem}[section]

\newtheorem{lem}[thm]{Lemma}
\newtheorem{prop}[thm]{Proposition}
 
\newtheorem{tet}{Teichm\"uller's Theorem} 
\newtheorem{sfmc}{Strebel's Frame Mapping Condition}

\theoremstyle{definition}
\newtheorem{defn}{Definition}[section]

\theoremstyle{remark}
\newtheorem{rem}{Remark}[section]

\newcommand{\diam}{\operatorname{diam}}
\newcommand{\dist}{\operatorname{dist}}

\newcommand{\cl}{\operatorname{cl}}
\renewcommand{\mod}{\operatorname{mod}}
\newcommand{\tl}{\tilde}
\newcommand{\wtl}{\widetilde}
\newcommand{\eps}{\epsilon}
\newcommand{\EE}{{\Ccal{E}}}
\newcommand{\tw}{{{\rm\boldsymbol T}}}
\newcommand{\Eps}{\operatorname{E}}

\newcommand{\aaa}{\rm}
\newcommand{\Cbb}[1]{{{\Bbb{#1}}}}
\newcommand{\Ccal}[1]{{{\cal{#1}}}}
\newcommand{\crit}{{{\rm\boldsymbol C}}}
\newcommand{\pr}{{\cal P}}
\newcommand{\hol}{{\rm{\boldsymbol H}}}
\newcommand{\bran}{{\rm{\boldsymbol B}}}

\renewcommand{\Re}{\operatorname{Re}}
\renewcommand{\Im}{\operatorname{Im}}

\numberwithin{equation}{section}
\newcommand{\thmref}[1]{Theorem~\ref{#1}}
\newcommand{\propref}[1]{Proposition~\ref{#1}}
\newcommand{\secref}[1]{\S\ref{#1}}
\newcommand{\lemref}[1]{Lemma~\ref{#1}}
 
\newcommand{\figref}[1]{Figure~\ref{#1}}

\def\SBIMSMark#1#2#3{
 \font\SBF=cmss10 at 10 true pt
 \font\SBI=cmssi10 at 10 true pt
 \setbox0=\hbox{\SBF Stony Brook IMS Preprint \##1}
 \setbox2=\hbox to \wd0{\hfil \SBI #2}
 \setbox4=\hbox to \wd0{\hfil \SBI #3}
 \setbox6=\hbox to \wd0{\hss
             \vbox{\hsize=\wd0 \parskip=0pt \baselineskip=10 true pt
                   \copy0 \break%
                   \copy2 \break%
                   \copy4 \break}}
 \dimen0=\ht6   \advance\dimen0 by \vsize \advance\dimen0 by 8 true pt
 \dimen2=\hsize \advance\dimen2 by .25 true in
 \ht6=0pt \dp6=0pt
 \setbox8=\vbox to \dimen0{\vfill \hbox to \dimen2{\hss \copy6}}
 \ht8=0pt \dp8=0pt \wd8=0pt
 \copy8
}

\begin{document}
\SBIMSMark{1998/5}{May 1998}{}

\addtolength{\evensidemargin}{-0.7in}
\addtolength{\oddsidemargin}{-0.7in}

\title[Attractor of  Renormalization]{The Attractor of  Renormalization and Rigidity of Towers of Critical Circle Maps}
\author{Michael Yampolsky}
\begin{abstract}
We demonstrate the existence of a global attractor $\cal A$ with a Cantor set structure 
for the renormalization of critical circle  mappings. The set $\cal A$ is invariant under  a
  generalized renormalization transformation, whose action on $\cal A$ is conjugate
to the two-sided shift. 
\end{abstract}
\maketitle


\section{Introduction}
The empirical discovery  of universality phenomena 
in the transition to chaos in one-dimensional dynamical systems during the
late 1970's  made a great impact  on the subject. To explain
these phenomena, parallels with statistical physics were drawn, and 
One-dimensional Renormalization Theory was born from this effort.
The main object of this theory is a renormalization transformation 
acting on an appropriate class of dynamical systems, and the universality
phenomena are related to the hyperbolic dynamics of the transformation.

The first renormalization transformation was constructed by Feigenbaum, and
Coullet \& Tresser in the setting of unimodal maps.
Unimodal renormalization theory enjoyed spectacular progress since
Sullivan introduced into it the methods of complex-analytic dynamics \cite{S1}.
The works of Sullivan \cite{Sul,MvS}, McMullen \cite{McM2}
and Lyubich \cite{Lyubich-feigenbaum,Lyubich-horseshoe} culminated in establishing the 
hyperbolicity of  the unimodal renormalization operator, thus providing the 
mathematical basis for the universality phenomena in unimodal dynamics.

The theory of renormalization of critical circle maps has developed alongside
with the unimodal theory. Its objects, the critical circle maps, are 
orientation preserving self-homeomorphisms of
 the circle ${\bf T}={\Bbb R}/{\Bbb Z}$ of class $C^3$ with a single critical
point $c$. A further assumption is made that the critical point is of cubic type.
This means that for a lift $\bar f:{\Bbb R}\to {\Bbb R}$ of a critical circle  map $f$
with critical points at integer translates of $\bar c$,
$$\bar f(x)-\bar f(\bar c)=(x-\bar c)^3(\operatorname{const}+O(x-\bar c)).$$
Examples of analytic critical circle maps are provided
by the projections $f_\theta$ to ${\Bbb R}/{\Bbb Z}$ of  homeomorphisms of the standard
 (or Arnold's) family 
$$A_\theta(x)=x+\theta-\frac{1}{2\pi}\sin 2\pi x.$$
 The rotation number of a critical circle map $f$ will be denoted by $\rho(f)$. 
The number theoretical properties of $\rho(f)$ have dynamical implications for $f$.
In particular, by a theorem
of Yoccoz \cite{Yoc}, in the case when $\rho(f)$ is irrational $f$ is topologically 
conjugate to the rigid rotation of the circle by angle $\rho(f)$.
It is useful to associate to a rotation number a continued fraction expansion 
$$
\rho(f)=\cfrac{1}{r_0+\cfrac{1}{r_1+\cfrac{1}{r_2+\dotsb}}}
$$
which we will further abbreviate as $\rho(f)=[r_0,r_1,\ldots]$ for typographical convenience.
The sequence $\{r_i\}$ is infinite if and only if $\rho(f)$ is irrational, in this
case it is said to be of {\it type bounded by $B$} if all terms $r_i$  are not greater than $B$.

The renormalization operator $\cal R$ for critical circle maps is defined in the 
language of critical commuting pairs. This approach first appeared in \cite{ORSS} and, in a slightly different form,
in \cite{FKS}, and has been further developed by Lanford \cite{La1,La2} and others.
Critical commuting pairs correspond to smooth conjugacy classes of critical circle
maps; in particular, they posess rotation numbers, on which renormalization acts
as the Gauss shift $G:\rho\to\{1/\rho\}$. All pairs with non-zero rotation numbers
are  renormalizable, those with irrational rotation numbers are infinitely renormalizable.
The main open question of the theory is the Hyperbolicity Conjecture, which
postulates the existence of a smooth structure on the space of commuting
pairs in which the renormalization transformation $\cal R$ is globally hyperbolic
with one-dimensional expanding direction. In its full generality the Conjecture is
due to Lanford (see \cite{La2}), below we discuss it in some more detail. 
Although the first attempts to prove the Conjecture were confined to the framework of one-dimensional
smooth dynamics, the importance of analytic methods was early understood.
In 1986 Eckmann and H. Epstein \cite{EE} constructed a class of real-analytic maps
invariant under $\cal R$, the so-called {\it Epstein class} $\cal E$. It was further
shown (see \cite{dF1}) that renormalizations of $C^3$ circle maps 
converge to $\cal E$, by a recent result of de~Faria and de~Melo \cite{dFdM1}
the convergence occurs at a geometric rate.

In \cite{dF1,dF2} de~Faria has adapted Sullivan's Renormalization Theory to the 
setting of critical circle mappings.
De~Faria has defined holomorphic extensions 
of renormalizations of maps in the Epstein class, which are the appropriate analogues 
of quadratic-like maps.
Using  Sullivan's techniques, he then demonstrated the existence of 
complex {\it a priori} bounds for such extensions in the case of maps of
bounded type, and obtained the  following renormalization convergence result:

\medskip
\noindent
{\bf Theorem \cite{dF1,dF2}.} {\sl Let $f_1$ and $f_2$ be two critical circle maps in $\cal E$
with the same rotation numbers in ${\Bbb R}\setminus{\Bbb Q}$ of bounded type. Then 
$$\operatorname{dist}_{C^r}({\cal R}^{\circ n}f_1,{\cal R}^{\circ n}f_2)\to 0 \text{ for all }0\leq r<\infty$$}

\medskip
As a consequence, for any $B\geq 1$ there exists a closed $\cal R$-invariant set ${\cal A}_B$
 of critical commuting pairs such that for any $f\in\cal E$ whose type is bounded by $B$,
${\cal R}^{\circ n}f\to {\cal A}_B$. In a recent work de~Faria and de~Melo \cite{dFdM2}
employed McMullen's towers techniques to  
demonstrate that the convergence to the attractor  ${\cal A}_B$ happens at a geometric rate.

In \cite{Ya1} we establish the existence of
complex {\it a priori} bounds for circle maps with an arbitrary irrational rotation
number, using methods developed in \cite{LY}. This enabled us to extend the above stated renormalization
convergence  result of de~Faria to circle maps with arbitrary irrational rotation numbers,
getting rid of the condition of bounded type.

In this paper we establish the existence of a global  attractor $\cal A$ for the 
renormalization operator $\cal R$, with a ``horseshoe" structure.
Let $\Sigma$ be the space of bi-infinite sequence of natural numbers, and 
denote by $\sigma:\Sigma\to\Sigma$ the shift on this space:
$$\sigma:(r_i)^\infty_{-\infty}\mapsto(r_{i+1})^\infty_{-\infty}.$$
For future use let us complement the natural numbers with the symbol $\infty$,
 and denote by $\bar \Sigma$ the space $({\Bbb N}\cup\{\infty\})^{\Bbb Z}$.
We prove the following:

\medskip
\noindent
{\bf Theorem A.} {\sl 
There exists a $\cal R$-invariant  set $\cal I$ of commuting pairs with irrational rotation numbers
with the following properties.
The action of $\cal R$ on $\cal I$ is bijective. Moreover,
there is a one-to-one correspondence
$$i:{\cal I}\to \Sigma$$
such that 
if $\zeta=i^{-1}(\ldots,r_{-k},\ldots,r_{-1},r_0,r_1,\ldots,r_k,\ldots)$ then
$\rho(\zeta)=[r_0,r_1,\ldots,r_k,\ldots]$ and thus the action of $\cal R$ on $\cal I$
is conjugate to the shift:
$$i\circ{\cal R}\circ i^{-1}=\sigma.$$
The set $\cal I$ is pre-compact in {\it Carath\'eodory topology} (see \secref{preliminaries} for 
the definiton), its closure $\cal A$ is the attractor for
the  renormalization operator:
$${\cal R}^{\circ n}\zeta\to {\cal A},\text{ for all }\;\zeta\in\cal E \text{ with }\rho(\zeta)\in{\Bbb R}\setminus{\Bbb Q},$$
where the convergence is understood in the sense of Carath\'eodory topology (it implies, in particular,
that the analytic extensions of the renormalized pairs converge uniformly on compact sets). More precisely,
for any pair $\zeta'\in \cal A$ with $\rho(\zeta)=\rho(\zeta')$ we have 
$$\dist ({\cal R}^{\circ n}\zeta,{\cal R}^{\circ n}\zeta')\to 0$$
for the $C^0$-distance between the analytic extensions of the renormalized pairs on an open neighborghood
of the origin.

}

\medskip
\noindent
The proof of the above theorem is inspired by the argument of McMullen for the convergence of unimodal
renormalizations \cite{McM2}. Following \cite{McM2} we consider the {\it geometric limits}
$T$ of various rescalings of sequences of renormalizations 
$(\zeta,{\cal R}\zeta,\ldots,{\cal R}^{\circ n}\zeta)$
of increasing lengths. The objects $T$   should be thought of as {\it bi-infinite towers} of 
nested dynamical systems in the plane. The proof of the theorem relies on the following uniqueness
statement:

\medskip
\noindent 
{\bf Tower Rigidity Theorem.}
{\sl There exists a unique (up to a homothety) bi-infinite tower 
for each bi-infinite sequence in $\bar \Sigma$.}

\medskip
\noindent
We point out that a parallel statement has been proved by  Hinkle \cite{Hinkle} in the setting of unimodal 
maps with essentially bounded combinatorics.
An important difference with the situation considered by McMullen and its adaptation to the case
of circle maps carried out in \cite{dFdM2} is the presense of {\it parabolic commuting pairs} in a limiting
tower. {\it Forward-infinite} towers with parabolic elements which occur as geometric limits in a broad
class of complex-analytic 
dynamical systems were considered in  the  dissertation of A. Epstein  \cite{Ep}, who proved
a rigidity theorem for such objects. It is conceivable that the constructions could be modified
appropriately so that the results of Epstein could be applied directly in the setting of critical circle maps;
we use some ergodic arguments to replace them.


The above mentioned parabolic elements of the attractor $\cal A$ are the pairs $\zeta$
in the closure of $\cal I$ 
with $\rho(\zeta)=0$ (we shall see below that such pairs posess fixed points with unit eigenvalues);
on these pairs the action of the renormalization operator $\cal R$ is not defined.
There is, however, a natural extension of ${\cal R}$ to commuting pairs 
with zero rotation number, the  {\it parabolic renormalization} $\cal P$. Combining the action of $\cal R$
and $\cal P$ into the {\it generalized renormalization operator } $\cal G$ we obtain:

\medskip
\noindent
{\bf Theorem B.} {\sl The attractor $\cal A$ is invariant under $\cal G$.
 The correspondence $i:{\cal I}\to \Sigma$ bijectively extends to 
$i:{\cal A}\to \bar\Sigma$ with the same properties, so that $i\circ {\cal G}\circ i^{-1}=\sigma$.}

\medskip
\noindent
There is, strictly speaking, no canonical way to extend  the parabolic renormalization operator
${\cal P}$ to Epstein commuting pairs with zero rotation number wich are not contained in the attractor $\cal A$.
The parabolic renormalization of a pair $\zeta$ is determined by the selection of  the rotation number
of ${\cal P}\zeta$. We may, however, arbitrarily associate a rotation number with every $\zeta\in\cal E$
having $\rho(\zeta)=0$ and a parabolic fixed point,
 and in this way extend the generalized renormalization operator $\cal G$ to the parabolic elements of 
$\cal E$. For any such extension we will have  the same convergence property:

\medskip
\noindent
{\bf Theorem C.} {\sl The generalized renormalizations ${\cal G}^n\zeta$ converge to $\cal A$ for
 any parabolic pair $\zeta\in\cal E$.}

\medskip
\noindent
As a consequence of the above theorems we obtain the 
following analogue of the golden mean renormalization fixed point:

\medskip
\noindent
{\bf Theorem D.} {\sl There exists a commuting pair $\zeta_0\in\cal A$ such that for all maps $f$ with 
$\rho(f)\in{\Bbb R}\setminus{\Bbb Q}$ and
$G^{\circ n}(\rho(f))\to 0$ we have
$${\cal R}^{\circ n}f\to \zeta_0.$$
Moreover, let ${\cal E}^0\subset \cal E$ be the set of Epstein pairs with zero rotation numbers,
and let ${\cal P}^0:{\cal E}^0\to{\cal E}^0$ denote the appropriate parabolic renormalization operator.
Then $\zeta_0$ is fixed under parabolic renormalization, ${\cal P}^0\zeta_0=\zeta_0$,
and 
$${\cal P}^0\zeta\to\zeta_0\text{ for all }\zeta\in{\cal E}^0.$$
}

\medskip
\noindent
{\bf Acknowledgements.} I would like to express my gratitude to A.~Epstein and B.~Hinkle,
the discussions with them have greatly assisted me in producing this paper. Ben Hinkle has 
a parallel work \cite{Hinkle} on unimodal maps with essentially bounded combinatorics,
and I appreciate his willingness to let me consult his manuscript before it was finished.
I wish to thank O. Lanford, whose numerous helpful comments on an earlier version of 
this paper helped me to streamline some of the arguments and to improve the exposition.
I also would like to thank E.~deFaria, who graciously provided me with preliminary
copies of his papers with W.~deMelo. Finally, my thanks go to Mikhail Lyubich for
numerous stimulating conversations and constant moral support.

\section{Preliminaries}
\label{preliminaries}

{\bf \noindent Some notations.}
The notation $D_r(z)$ will stand for the Euclidean disk with center $z\in\Bbb C$ and
radius $r$. The unit disk $D_1(0)$ will be denoted $\Bbb D$. 
For two points $a$ and $b$ in the circle $\bf T$, $[a,b]$ will denote the shorter
of the two arcs connecting them; $|[a,b]|$ will denote the length of the 
arc. For two points $a,b\in \Bbb R$, $[a,b]$ will denote the closed interval with
endpoints $a$, $b$ without specifying their order.
The plane $({\Bbb C}\setminus{\Bbb R})\cup J$
with the parts of the real axis not contained in the interval $J\subset \Bbb R$
removed  will be denoted ${\Bbb C}_J$.

We use $\dist$ and $\diam$ to denote the Euclidean distance and diameter in $\Bbb C$.
We call two real numbers $x$ and $y$ {\it $K$-commensurable} or simply {\it commensurable}
if $K^{-1}\leq |x|/|y|\leq K$ for some $K>1$. Two sets $X$ and $Y$ in $\Bbb C$ are 
$K$-commensurable if their diameters are.
In accordance with the established terminology, we shall say that a quantity is
{\it definite} if it is greater than a universal positive constant.
 A set $B$ is contained {\it well inside } of a set
$A\subset \Bbb C$ if  $ A\setminus B$ contains an annulus with definite modulus.
Similarly, an interval $I\subset \Bbb R$ is contained well inside of another interval $J$
there exists a universal constant $K>0$ such that for each component $L$ of $J\setminus I$
we have $|L|>K|I|$.

\medskip

{\bf \noindent Commuting pairs and renormalization of critical circle maps.}
Consider a critical circle mapping $f$ with a rotation number $\rho$, and let
\begin{equation}
\label{rotation-number}
\rho(f)=[r_0,r_1,r_2,\ldots]
\end{equation}
be its (possibly finite) continued fraction expansion. 
We shall always assume that the critical point of $f$ is at $0\in{\Bbb R}/{\Bbb Z}$.
An iterate $f^k(0)$ is called a {\it closest return} of the critical point if
the arc $[0,f^k(0)]$ contains no other iterates $f^i(0)$ with $i<k$.
Since  a circle homeomorphism $f$ is semi-conjugated
to the rigid rotation by angle $\rho(f)$ (see for example \cite{MvS}), 
the moments of closest returns are determined by the number-theoretic properties of $\rho(f)$.
Namely, the closest returns occur at iterates $\{f^{q_m}(0)\}$ where $q_m$'s are given recursively
 by $q_{m+1}=r_mq_m+q_{m-1}$, $q_0=1$, $q_1=r_0$. Thus the
number $q_m$ appears as the denominator of the  truncated continued
fraction expansion of $\rho$ of length $m-1$ in its reduced form:
\begin{equation}
\label{rational-approximation}
p_m/q_m=[r_0,r_1,\ldots,r_{m-1}].
\end{equation}
Set $I_m\equiv[0,f^{q_m}(0)]$. As a consequence of
 \'Swia\c\negthinspace tek-Herman real {\it a priori} bounds
(\cite{Sw,H}), the intervals $I_m$ and $I_{m+1}$ are $K-$commensurable, 
with a universal constant $K$ provided $m$ is large enough.

The general strategy of defining a renormalization of a given dynamical system 
(we are following \cite{Lyubich-cambridge} here) is to select a piece of its phase
space, rescale it to the ``original'' size, and then consider the first return map
to this piece. Historically, for a circle map $f$ the union of arcs $A_m=I_m\cup I_{m+1}$
is chosen as the domain for the return map. The first return map $R_m:A_m\to A_m$ is 
defined piecewise by by $f^{q_m}$ on $I_{m+1}$  and by $f^{q_{m+1}}$ on $I_m$.
To view $R_m$ as a critical circle map we may identify the neighborhoods of points 
$f^{q_m}(0)$ and $f^{q_{m+1}}(0)$ by the iterate $f^{q_{m+1}-q_m}$. This identification
transforms the arc $A_m$ into  a $C^3$-smooth closed one-dimensional manifold $\tl A_m$,
$R_m$ projects to a smooth homeomorphism $\tl R_m:\tl A_m\to \tl A_m$ with a critical point at $0$.
However, the manifold $\tl A_m$ does not posess a canonical affine structure;
the choice of a smooth identification $\phi:\tl A_m\to {\Bbb R}/{\Bbb Z}$ gives rise to 
a plethora of different critical circle maps, all smoothly conjugate.
The above discussion illustrates why the space of critical circle maps is not 
suited for defining a renormalization transformation, and motivates the introduction of 
the following objects:
 
\begin{defn}
A  {\it  commuting pair} $\zeta=(\eta,\xi)$ consists of two 
smooth  orientation preserving interval homeomorphisms 
$\eta:I_\eta\to \eta(I_\eta),\;
\xi:I_{\xi}\to \xi(I_\xi)$, where
\begin{itemize}
\item[(I)]{$I_\eta=[0,\xi(0)],\; I_\xi=[\eta(0),0]$;}
\item[(II)]{Both $\eta$ and $\xi$ have homeomorphic extensions to interval
neighborhoods of their respective domains with the same degree of
smoothness, which commute, 
$\eta\circ\xi=\xi\circ\eta$;}
\item[(III)]{$\xi\circ\eta(0)\in I_\eta$;}
\item[(IV)]{$\eta'(x)\ne 0\ne \xi'(y) $, for all $x\in I_\eta\setminus\{0\}$,
 and all $y\in I_\xi\setminus\{0\}$.}
\end{itemize}
\end{defn}

\noindent
Commuting pairs were first used to define the renormalization transformation
by Ostlund, Rand, Sethna, and Siggia \cite{ORSS}.
 Feigenbaum, Kadanoff, and Shenker \cite{FKS} defined renormalization by means of
a slightly different formalism.

A {\it critical commuting pair} is a commuting pair $(\eta,\xi)$ whose maps  can be decomposed near zero as
$\eta=h_\eta\circ Q\circ H_\eta$,
and $\xi=h_\xi\circ Q\circ H_\xi$, where $h_\eta,h_\xi,
H_\eta,H_\xi$ are real analytic  diffeomorphisms and $Q(x)=x^3$. We shall further require 
a technical assumption that
$\xi$  analytically extends to an interval $(a,b)\ni 0$ with $\xi(a,b)\supset [\eta(0),\xi(0)]$,
and has a single critical point $0$ in this interval.
The space of critical commuting pairs modulo affine conjugacy endowed with 
$C^0$ topology, will be denoted by $\crit$.

Let $f$ be a critical circle mapping, whose rotation number $\rho$
has a continued fraction expansion (\ref{rotation-number}) with
at least $m+1$ terms. Let $p_m$ and $q_m$ be as in 
(\ref{rational-approximation}). The pair of iterates $f^{q_{m+1}}$
and $f^{q_m}$ restricted to the circle arcs $I_m$ and $I_{m+1}$
correspondingly can be viewed as a critical commuting pair  in
the following way.
Let $\bar f$ be the lift of $f$ to the real line satisfying $\bar f '(0)=0$,
and $0<\bar f (0)<1$. For each $m>0$ let $\bar I_m\subset \Bbb R$ 
denote the closed 
interval adjacent to zero which projects down to the interval $I_m$.
Let $\tau :\Bbb R\to \Bbb R$ denote the translation $x\mapsto x+1$.
Let $\eta :\bar I_m\to \Bbb R$, $\xi:\bar I_{m+1}\to \Bbb R$ be given by
$\eta\equiv \tau^{-p_{m+1}}\circ\bar f^{q_{m+1}}$,
$\xi\equiv \tau^{-p_m}\circ\bar f^{q_m}$. Then the pair of maps
$(\eta|\bar I_m,\xi|\bar I_{m+1})$ forms a critical commuting pair
corresponding to $(f^{q_{m+1}}|I_m,f^{q_m}|I_{m+1})$.
Henceforth we shall  simply denote this commuting pair by
\begin{equation}
\label{real1}
(f^{q_{m+1}}|I_m,f^{q_m}|I_{m+1}).
\end{equation}
This allows us to readily identify the dynamics of
the above commuting pair with that of the underlying circle map,
at the cost of a minor abuse of notation.

Following \cite{dFdM1}, we say that the {\it height} $\chi(\zeta)$
of a critical commuting pair $\zeta=(\eta,\xi)$ is equal to $r$,
if 
$$0\in [\eta^r(\xi(0)),\eta^{r+1}(\xi(0))].$$
 If no such $r$ exists,
we set $\chi(\zeta)=\infty$, in this case the map $\eta|I_\eta$ has a 
fixed point.  For a pair $\zeta$ with $\chi(\zeta)=r<\infty$ one verifies directly that the
mappings $\eta|[0,\eta^r(\xi(0))]$ and $\eta^r\circ\xi|I_\xi$
again form a commuting pair.
For a commuting pair $\zeta=(\eta,\xi)$ we will denote by 
$\wtl\zeta$ the pair $(\wtl\eta|\wtl{I_\eta},\wtl\xi|\wtl{I_\xi})$
where tilde  means rescaling by a linear factor $\lambda={1\over |I_\eta|}$.

\begin{defn}
The {\it renormalization} of a real commuting pair $\zeta=(\eta,
\xi)$ is the commuting pair
\begin{center}
${\cal{R}}\zeta=(
\widetilde{\eta^r\circ\xi}|
 \widetilde{I_{\xi}},\; \widetilde\eta|\widetilde{[0,\eta^r(\xi(0))]}).$
\end{center}
\end{defn}

\noindent
The non-rescaled pair $(\eta^r\circ\xi|I_\xi,\eta|[0,\eta^r(\xi(0))])$ will be referred to as the 
{\it pre-renormalization} $p{\cal R}\zeta$ of the commuting pair $\zeta=(\eta,\xi)$.

For a pair $\zeta$ we define its {\it rotation number} $\rho(\zeta)\in[0,1]$ to be equal to the 
continued fraction $[r_0,r_1,\ldots]$ where $r_i=\chi({\cal R}^i\zeta)$. 
In this definition $1/\infty$ is understood as $0$, hence a rotation number is rational
if and only if only finitely many renormalizations of $\zeta$ are defined;
if $\chi(\zeta)=\infty$, $\rho(\zeta)=0$.
Thus defined, the rotation number of a commuting pair can be viewed as a rotation number in
the usual sense via the following construction.
Given a critical commuting pair $\zeta=(\eta,\xi)$
we can regard the interval $I=[\eta(0),\xi\circ\eta(0)]$ as a circle, identifying 
$\eta(0)$ and $\xi\circ\eta(0)$ and define $f_\zeta:I\to I$ by 
$$f_\zeta=\left\{\begin{array}{l}
                    \xi\circ\eta(x) \text{ for }x\in [\eta(0),0]\\
                    \eta(x)\text{ for } x\in [0,\xi\circ\eta(0)]
                \end{array}\right.
$$
We perform glueing together of $\eta(0)$ to $\xi\circ\eta(0)$ by the mapping
$\xi$, which by the condition (II) above extends to a smooth homeomorphism
of open neighborhoods. The quotient of the interval $I$ is a 
closed one-dimensional manifold $M$, the mapping $f_\zeta$ projects down
to a smooth homeomorphism $F_\zeta:M\to M$.
Identifying $M$ with the circle by a diffeomorphism $\phi:M\to S^1$
we recover a critical circle mapping $f^\phi=\phi\circ F_\zeta\circ\phi^{-1}$.
The critical circle mappings corresponding to two different choices
of $\phi$ are conjugated by a diffeomorphism, and thus we recovered
a smooth conjugacy class of circle mappings from a critical commuting
pair.
It is immediately seen, that:
\begin{prop}
\label{rotation number}
The rotation number of mappings in the above constructed conjugacy class
is equal to $\rho(\zeta)$.
\end{prop}

\noindent
The advantage of defining $\rho(\zeta)$ using a sequence of heights is
in having a way of distinguishing the commuting pairs with rotation numbers $0$ and $1$,
this way we also remove the ambiguity in 
prescribing a continued fraction expansion to rational rotation numbers.

For $\rho=[r_0,r_1,\ldots]\in [0,1]$ let us set 
$$G(\rho)=[r_1,r_2,\ldots]=\left\{\frac{1}{\rho}\right\},$$
where $\{x\}$ denotes the fractional part of a real number $x$
($G$ is usually referred to as the {\it Gauss map}).
As follows from the definition,
$$\rho({\Ccal R}\zeta)=G(\rho(\zeta))$$ for a real commuting pair $\zeta$ with
$\rho(\zeta)\ne 0$.

The renormalization of the real commuting pair (\ref{real1}), associated
to some critical circle map $f$, is the rescaled pair
$(\wtl{{f}^{q_{m+2}}}|\wtl{{I}_{m+1}},\wtl{{f}^{q_{m+1}}}|\wtl{{I}_{m+2}})$.
Thus for a given critical circle map $f$ the renormalization operator
 recovers the (rescaled) sequence of the first return maps:
$$\{ (\wtl{f^{q_{i+1}}}|\wtl{I_i},\wtl{f^{q_{i}}}|\wtl{I_{i+1}})\}_{i=1}^{\infty}.$$

\comm{
Let $\zeta=(\eta,\xi)$ be a commuting pair with rotation number $0$,
in which case the interval map $\eta|_{I_\eta}$ has a fixed point.
We will distinguish the cases $\eta(0)>0$ and $\eta(0)<0$, by writing 
$\rho(\zeta)=0_-$ and $\rho(\zeta)=0_+$ correspondingly, referring
to $\rho$ in this case as the {\it signed rotation number of $\zeta$}.
 For a pair $\zeta$ with a rational $\rho$, 
$\rho=p/q=[r_0,r_1,\ldots,r_n]$ we will write $\rho={p/q}_-$ when 
$\rho({{\cal R}^{n+1}}\zeta)=0_-$ and similarly for ${p/q}_+$.
}

Let us denote by $\crit_\infty\subset \crit$ the space of critical commuting pairs $\zeta$
modulo rescaling
with $\chi(\zeta)=\infty$. It is clear that ${\cal R}:\crit\setminus\crit_\infty\to\crit$.
Let us note:

\begin{rem}
\label{injective1}
The map ${\cal R}:\crit\setminus\crit_\infty\to\crit$ is injective.
\end{rem}
\begin{pf}
Let $\zeta=(\eta,\xi)$ be a pair in ${\cal R}(\crit\setminus\crit_\infty)$. 
For $r\in \Bbb N$ let $\gamma_r$ be the maximal analytic extension of
$\xi^{-r}\circ\eta$. Choose the smallest $r_0$ for which $\gamma_{r_0}^{-1}$
is differentiable at $\xi(0)$. 
 Setting $\gamma=\gamma_{r_0}$ and
$$\zeta_{-1}=(\xi|_{[0,\gamma(0)]},\gamma|_{[0,\xi(0)]})$$
we have $\zeta_{-1}\in\crit$ and ${\cal R}\zeta_{-1}=\zeta$.
\end{pf}

\medskip
\noindent
{\bf Renormalization hyperbolicity conjecture.}
We would like to discuss briefly the general Renormalization Hyperbolicity conjecture 
and its connection to the existence of the attractor of the renormalization operator.

\realfig{cylinder}{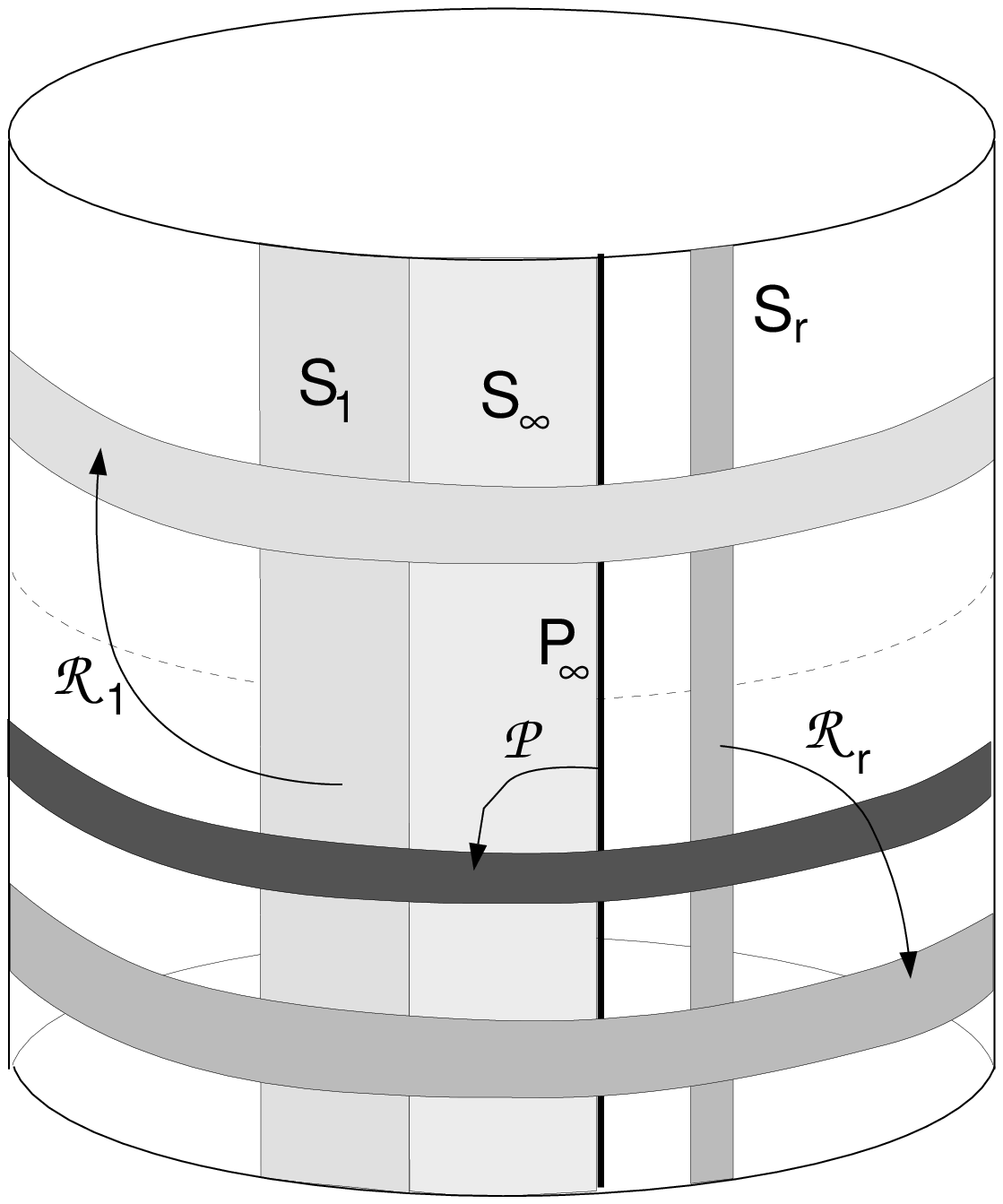}{}{6.1cm}

\smallskip
\noindent
{\bf Renormalization Hyperbolicity Conjecture.} {\sl There exists a renor\-ma\-liza\-tion-invariant
set of critical commuting pairs with the structure of an infinite dimensional smooth
manifold, 
with respect to which the renormalization transformation $\cal R$ is globally
uniformly hyperbolic, with one-dimensional expanding direction.

Moreover, if $t$ denotes the expanded coordinate in a local chart, then the dependence
$t\mapsto \rho(\zeta(t,\cdot))$ is continuous and (not strictly) monotone.}

\smallskip
\noindent
In this generality the conjecture is due to Lanford \cite{La2}, it remains the main open question
of the theory.
To illustrate the conjecture it is convenient to consider the following caricature.
Let $\rho=\phi(\theta):{\bf T}\to{\bf T}$ be a monotone continuous function such that
$\phi^{-1}(\rho)$ is a point for $\rho$ irrational and an interval otherwise.
An example of such a function is the dependence $\theta\mapsto\rho(f_\theta)$ of the 
rotation number of a standard map on the parameter.
Imagine the relevant space of commuting pairs as an infinite-dimensional cylinder
${\aaa C}={\bf T}\times{\aaa C}'$, where the rotation number of a commuting
pair $\zeta(\theta,\cdot)$ with the equatorial coordinate $\theta\in{\bf T}$ is $\rho(\theta)$.
 The cylinder $\aaa C$ is partitioned into strips  (cf. \figref{cylinder})
$${\aaa S}_r=\{\zeta\in{\aaa C}|\rho(\zeta)=[r,r_1,\ldots]\}\text{ for }r=1,\ldots,\infty$$
A boundary component of the strip ${\aaa S}_\infty$ is the hypersurface ${\aaa P}_\infty\subset 
{\aaa S}_\infty$ with the property that a pair $\zeta\in {\aaa P}_\infty$ if and only if
it has a fixed point with unit eigenvalue.
The sets ${\aaa S}_r$ accumulate on ${\aaa P}_\infty$ in clockwise direction.

It is natural to think of the transformation $\cal R:{\aaa C}\setminus {\aaa S}_\infty\to {\aaa C}$ 
as being defined piecewise, given on
each ${\aaa S}_r$, $r\ne \infty$ by the formula:
$${\cal R}_r:(\eta,\xi)\mapsto(\eta,\eta^r\circ\xi).$$
The operator ${\cal R}_r$ uniformly expands  the strip ${\aaa S}_r$ in the equatorial direction,
and uniformly contracts it in all other directions, mapping it  onto a thin cylinder intersecting all the strips.
The invariant set $\cal I$ is seen in this picture as intersections of the ``boxes''
$$
{\cal R}_{r_{-1}}({\aaa S}_{r_{-1}}\cap {\cal R}_{r_{-2}}({\aaa S}_{r_{-2}}\cap\cdots ({\cal R}_{r_{-n}}{\aaa S}_{r_{-n}})
\cdots)\cap {\aaa S}_{r_0}\cap {\cal R}^{-1}_{r_0}({\aaa S}_{r_1}\cap {\cal R}_{r_1}^{-1}({\aaa S}_{r_2}\cap \cdots
{\cal R}_{r_{n-1}}^{-1}({\aaa S}_{r_n})\cdots)
$$
The parabolic renormalization ${\cal P}$ which we define below transforms the set ${\cal P}_\infty$ into a 
thin equatorial cylinder, the hyperbolicity conjecture can be extended in the obvious way to include
this transformation.

\medskip

\noindent
{\bf Holomorphic commuting pairs.}
\realfig{circ3}{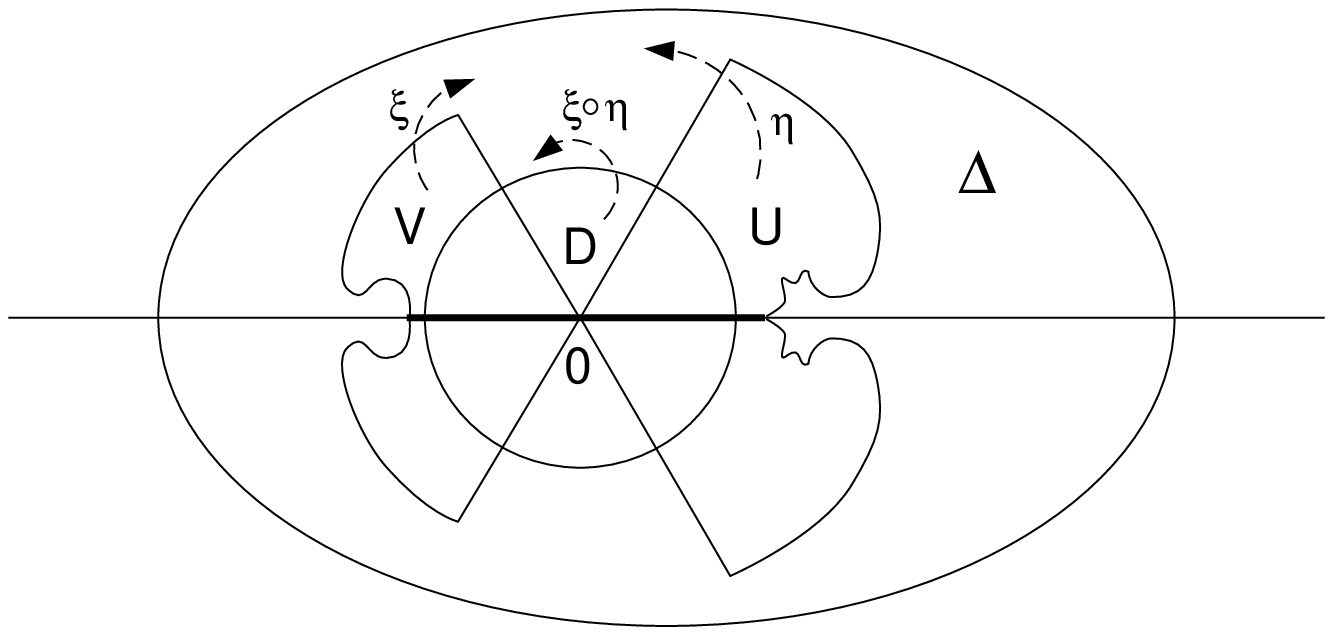}{}{10cm}
Following \cite{dF1,dF2} we say that a real commuting pair
$ ( \eta , \xi )$
extends to a {\it holomorphic commuting pair}  $\cal H$
 (cf. \figref{circ3})
if there exist four 
$\Bbb R$-symmetric domains $\Delta$, $D$, $U$, $V$, and a holomorphic 
mapping $\nu$, such that 
\begin{itemize}
\item  $\bar D,\; \bar U,\; \bar V\subset \Delta$,
 $\bar U\cap \bar V=\{ 0\}$;
 $U\setminus D$,  $V\setminus D$, $D\setminus U$, and $D\setminus V$ 
 are nonempty 
connected sets, $U\supset I_\eta$, $V\supset I_{\xi}$;
\item mappings $\eta:U\to \Delta\cap {\Bbb C}_{\eta(J_{U})}$ and
 $\xi:V\to \Delta\cap{ \Bbb C}_{{\xi}(J_{V})}$ are onto and
univalent, where $J_{U}=U\cap \Bbb R$, $J_{V}=V\cap \Bbb R$;
\item $\nu:D\to \Delta\cap{\Bbb C}_{\nu(J_D)}$ is a three-fold 
branched covering with a unique fixed point at zero, where 
$J_D=D\cap {\Bbb R}$;
\item $\eta$ and $\xi$ have holomorphic extensions to
a certain neighborhood of the origin where both $\eta\circ \xi$
and $\xi\circ\eta$ are defined, and $\eta\circ \xi(z)=\xi\circ\eta(z)=\nu(z)$
for all $z$ in that neighborhood.
\end{itemize}

As the  following proposition shows, the mapping $\nu$ is 
nothing else, but the composition of $\eta$ and $\xi$:
\begin{prop}[Proposition II.1,\cite{dF2}]
Under the above conditions, the mappings $\eta$ and $\xi$ have analytic 
extensions to $U\cup D$ and $V\cup D$ correspondingly.
Moreover, $\eta:D\to V\cap {\Bbb C}_{[\xi^{-1}(\eta(0)),0]}$
and $\xi:D\to U$ are three-fold branched covering maps,
and $\nu=\eta\circ\xi$.
\end{prop}

Set $\Omega=D\cup U\cup V$.
One immediately observes that if a commuting pair $\zeta=(\eta,\xi)$ 
with a finite height has
a holomorphic pair extension ${\cal H}:\Omega\to \Delta$, then there exists a 
holomorphic commuting pair ${\cal G}:\Omega'\to \Delta$  whose restriction
coincides with ${\cal R}\zeta$. We shall refer to the commuting pair $\zeta$
as {\it underlying} $\cal H$, and write $\zeta={\cal H}\cap {\Bbb R}$.

We say that a real commuting pair $(\eta,\xi)$ with
an irrational rotation number has
{\it complex bounds}, if all its renormalizations extend to 
holomorphic commuting pairs with definite moduli, that is
$$\mod(\Delta\setminus\Omega)>\mu>0.$$

 The {\it shadow} of the holomorphic commuting
pair is the piecewise holomorphic mapping $S_{\cal H}:\Omega\to \Delta$, 
given by
$$S_{\cal H}(z) = \left\{ \begin{array}{l}
		   \eta(z),\; z\in U\\
		   \xi(z),\; z\in V\\
		   \xi\circ\eta(z),\; z\in D\setminus (U\cup V)
	         \end{array} \right. $$
The shadow of a holomorphic pair captures its dynamics in the following sense:
\begin{prop}[Prop. II.4. \cite{dF2}]
\label{shadow-cover}
Given a holomorphic commuting pair $\cal H$ as above, consider its shadow 
$S_{\cal H}$.
Let $I=\Omega \cap \Bbb R$, and $X=I\cup S_{\cal H}^{-1}(I)$. Then:
\begin{itemize}
\item The restriction of $S_{\cal H}$ to $\Omega\setminus X$ 
is a regular three fold
covering onto $\Delta\setminus \Bbb R$.
\item The orbits of $S_{\cal H}$ and $\cal H$ coincide as sets.
\end{itemize}
\end{prop}

We will say that two holomorphic commuting pairs
${\cal H}:\Omega_{\cal H}\to D_{\cal H}$ and 
${\cal G}:\Omega_{\cal G}\to D_{\cal G}$ are conjugate if there 
is a homeomorphism $h:D_{\cal G}\to D_{\cal H}$ such that
$$S_{\cal G}=h^{-1}\circ S_{\cal H}\circ h.$$
We will usually write simply ${\cal G}=h^{-1}\circ {\cal H}\circ h,$
meaning that $h$ conjugates the corresponding elements of the two holomorphic pairs.
We define the filled Julia set $K(\cal H)$ of a holomorphic commuting
pair $\cal H$ as the closure of the collection of all points which do not escape
$\Omega$ under iteration of $S_{\cal H}$. This set is clearly compact and
connected, it is full by the Maximum principle. Its boundary is the 
Julia set of $\cal H$, denoted by $J({\cal H})$.

\medskip
\noindent
{\bf Holomorphic commuting pairs in the standard family.}
For each $0\leq \theta<1$ let $A_\theta$ be the entire mapping given by
$$A_\theta(z)=z+\theta+\frac{1}{2\pi}\sin(2\pi z).$$
Since $A_\theta\circ T=T\circ A_\theta$, where $T$ is the unit translation
$z\mapsto z+1$, each $A_\theta$ is a lift of a holomorphic self-mapping
of the cylinder, $f_\theta:\Bbb C/\Bbb Z\cong \Bbb C^*$.
As each $A_\theta$ is real-analytic and satisfies $A'_\theta(x)>0$
for $x\in \Bbb R\setminus \Bbb Z$, the restriction $f_\theta|\bf T$
is  a critical circle map, whose rotation number we will denote 
$\rho(\theta)$. These restrictions are  usually referred to as 
the standard family (or Arnold family) of circle homeomorphisms.
Elementary considerations of monotone dependence on parameter imply that 
 $\theta\to\rho(\theta)$ is a
continuous non-decreasing map of $\bf T$ onto itself. Whenever
$t\in \bf T$ is irrational, $\rho^{-1}(t)$ is a single point. For $t=p/q$ the set 
$\rho^{-1}(t)$ is a closed interval, for every parameter value in this interval
the homeomorphism $f_\theta$ has a period $q$ orbit. This orbit has eigenvalue
one at the two endpoints of the interval.

De Faria (cf. \cite{dF1,dF2}) demonstrated, using some explicit estimates 
on the growth of the standard maps, that renormalizations of maps
$f_\theta$ can be extended to holomorphic commuting pairs.
Before formulating his result, let us
 briefly describe the geometry of a mapping $A_\theta$.
The preimage of the real axis under $A_\theta$ consists of the 
axis itself together with the family of analytic curves
$$\Gamma_\pm^k:\Re(z)=k\pm \frac{1}{2\pi}\arccos\frac
{-2\pi|\Im(z)|}{\sinh(2\pi \Im(z))},$$ 
where $k\in \Bbb Z$. For each $k$ the curves $\Gamma^k_\pm$ meet
at the critical point $c_k=k$ and are both asymptotic to the vertical
lines $\Re(z)=k\pm 1/4$. Note that each $c_k$ is of cubic type.
The curves $\Gamma^k_\pm$ and $\Bbb R$ partition the complex plane
into simply-connected regions each of which is univalently mapped
onto $\Bbb H$ or $-\Bbb H$ by $A_\theta$.
Now denote by $U_n$ the connected component of the preimage 
$(A_\theta^{q_n})^{-1}(\Bbb H)$ whose boundary  contains the point 
$T^{-p_{n+1}}(A_\theta^{q_{n+1}}(0))$.
The closure of the union of $U_n$ and its reflection in $x-$axis
will be denoted $\hat U_n$. Similarly let $V_n$ be the component of
$(A_\theta^{q_{n+1}})^{-1}(\Bbb H)$ with 
$\cl V_n\ni T^{-p_{n}}(A_\theta^{q_{n}}(0))$, and set $\hat V_n$ to
be the union of $\cl V_n$ with its vertical reflection.

\begin{lem}[\cite{dF1,dF2}] 
\label{standard-pairs}
Set $\eta=T^{-p_{n+1}}\circ A_\theta^{q_{n+1}}$
and $\xi=T^{-p_{n}}\circ A_\theta^{q_{n}}$. For all sufficiently large
$R$ the preimages $U_{n,R}=\eta^{-1}(D_R(0))\cap \hat U_n$ and
$V_{n,R}=\xi^{-1}((D_R(0))\cap \hat V_n)$ are compactly contained in
$D_R(0)$. Thus the pair of maps $(\eta,\xi)$ extends to  a holomorphic pair with
range $D_R(0)$. Moreover the modulus of this holomorphic pair tends to infinity
with $R$.
\end{lem}

The significance of these examples in the Renormalization Theory is due to a 
rigidity property of standard maps:
\begin{lem}[see Lemma IV.8 \cite{dF2}]
\label{topological completeness}
Every real-symmetric holomorphic self-map of ${\Bbb C}/\Bbb Z$ which is topologically 
conjugate to a member of the family $\{f_\theta\}$ belongs to this family itself.
\end{lem}

\noindent
Note, that  $f_\theta$ is a finite-type analytic map of the cylinder, and
thus, by a  a theorem of L. Keen \cite{Keen}, does not have  wandering domains.
Combining \lemref{topological completeness} with the above described properties
of the dependence $\theta\mapsto\rho(\theta)$ we have the following:
\begin{lem}
\label{rigidity of standard maps}
Let $g$ be a real-symmetric self map of the cylinder ${\Bbb C}/\Bbb Z$ topologically
conjugate to a map $f_\theta$. If $\rho(\theta)$ is irrational, or 
$f_\theta$ has a periodic orbit with eigenvalue one on the circle, then
$g\equiv f_\theta$.
\end{lem}

\noindent
Taking further the first statement of  the above Lemma, de~Faria has shown:
\begin{lem}[\cite{dF1,dF2}]
\label{rigidity-irrational}
If $\rho(\theta)$ is irrational, then $f_\theta$ admits no non-trivial,
symmetric, invariant Beltrami differentials entirely supported on
its Julia set.
\end{lem}

\noindent
The above properties of standard maps carry over to general holomorphic commuting
pairs via quasiconformal {\it straightening} arguments.

\comm{
In the case when $\rho(\theta)$ is rational, let us single out the 
possibility when $f_\theta$ has a parabolic periodic point. In this
case, the filled Julia set of $f_\theta$ has interior, which
may support invariant Beltrami differentials, however,
\begin{lem}[\cite{EKT}]
\label{rigidity-parab}
If $f_\theta$ has a parabolic periodic point, there exist no non-trivial,
symmetric
quasiconformal deformations of $f_\theta$, entirely supported on its filled
Julia set.
\end{lem}

This second Lemma is a direct corollary of the following
\begin{thm}[\cite{EKT}]
\label{unique-parab}
For each signed rational number $p/q_\pm\in[0,1)$ there exists a 
unique $\theta\in[0,1)$ such that $f_\theta$ is a paraboilc circle map with 
signed rotation number $p/q_\pm$.
\end{thm}
The above Theorem is based on an ingenious Teichm\"uller theory argument,
originating in \cite{Ep},
which will be employed  further
 in this paper to prove a more general statement.
}

\medskip

\noindent
{\bf A topology on a space of branched coverings.}
Consider the  collection $\bran$ of all triplets $(U,u,f)$, where $U\subset \Cbb C$ is a topological
disk different from the whole plane, $u\in U$, and $f:U\to \Cbb C$ is a three-fold analytic 
branched covering map, with the only branch point at $u$.
We will endow $\bran$ with a topology as follows (cf. \cite{McM1}).

Let $\{(U_n,u_n)\}$ be a sequence of open connected regions 
$U_n\subset \Cbb C$ with {\it marked points} $u_n\in U_n$.
Recall that this sequence {\it Carath\'eodory
converges} to a marked region $(U,u)$ if:
\begin{itemize}
\item $u_n\to u\in  U$, and
\item for any Hausdorff limit point $K$ of the sequence 
$\hat{\Cbb C}\setminus U_n$, $U$ is a component of 
$\hat{\Cbb C}\setminus K$.
\end{itemize}
For a simply connected $U\subset \Cbb C$ and $u\in U$ let $R_{(U,u)}:\Cbb D\to U$ denote
the Riemann mapping with normalization $R_{(U,u)}(0)=u$, $R'_{(U,u)}(0)>0$.
By a classical result of Carath\'eodory, the Carath\'edory convergence
of simply-connected regions
$(U_n,u_n)\to (U,u)$ is equivalent to the locally uniform convergence of 
Riemann mappings $R_{(U_n,u_n)}$  to $R_{(U,u)}$.

For positive numbers $\eps_1$, $\eps_2$, $\eps_3$ and compact subsets $K_1$ and $K_2$
of the open unit disk $\Cbb D$, let the neighborhood 
$\Ccal U_{\eps_1,\eps_2,\eps_3,K_1,K_2}(U,u,f)$ of an element $(U,u,f)\in \bran$
be the set of all $(V,v,g)\in\bran$, for which:
\begin{itemize}
\item $|u-v|<\eps_1$,
\item $\displaystyle \sup_{z\in K_1}|R_{(V,v)}(z)-R_{(U,u)}(z)|<\eps_2$,
\item and $R_{(U,u)}(K_2)\subset V$, and $\displaystyle \sup_{z\in R_{(U,u)}(K_2)}|f(z)-g(z)|<\eps_3$.
\end{itemize}
One verifies that the sets $\Ccal U_{\eps_1,\eps_2,\eps_3,K_1,K_2}(U,u,f)$ form
a base of a topology on $\bran$, which we will call {\it Carath\'eodory topology}. 
This topology is clearly Hausdorff, and  the convergence of a sequence $(U_n,u_n,f_n)$ to $(U,u,f)$
is equivalent Carath\'eodory convergence of the marked regions $(U_n,u_n)\to (U,u)$ as well
as locally uniform convergence $f_n\to f$.

\medskip

\noindent
{\bf Epstein class.}
An orientation preserving interval homeomorphism $g:I=[0,a]\to g(I)=J$ 
belongs to the {\it Epstein class $\Ccal E$} if it extends to an
analytic three-fold branched covering map of a topological disk
 $G\supset I$ 
onto the double-slit plane ${\Bbb C}_{\tl J}$, where $\tl J\supset \cl J$.
Any map $g$ in the Epstein class can be decomposed as
\begin{equation}
\label{epstein-decomposition}
g=Q_c\circ h,
\end{equation}
where $Q_c(z)=z^3+c$, and $h:I\to [0,b]$ is a univalent map 
$h:G\to \Delta(h)$ onto the complex plane with six slits,
which triple covers ${\Bbb C}_{\tl J}$ under the cubic map $Q_c(z)$.

For any $s\in (0,1)$, let us introduce a smaller class
 ${\Ccal E}_s\subset \Ccal E$ of Epstein mappings
$g:I=[0,a]\to J\subset \tl J$
for which both $|I|$ and $\dist(I,J)$ are $s^{-1}$-commensurable with $|J|$,
the length of each component of $\tl J\setminus J$
is at least $s|J|$, and $g'(a)>s$.
We will often refer to the space $\Ccal E$ as {\it the} Epstein class,
and to each ${\Ccal E}_s$ as {\it an} Epstein class.

We say that a commuting pair $(\eta,\xi)\in\crit$
belongs to the (an) Epstein class if both of its maps do.
It immediately follows from the definitions that:
\begin{lem} 
The space of commuting pairs in the Epstein class ${\Ccal E}$ 
is invariant under renormalization.
\end{lem}
We shall denote by $\hol$ the space of holomorphic commuting pairs
$\Ccal H:\Omega\to \Delta$ whose underlying real commuting pair $(\eta,\xi)$
is in the Epstein class. 
In this case both maps $\eta$ and $\xi$ extend to triple branched coverings 
$\hat{\eta}:\hat{U}\to \Delta\cap {\Bbb C}_{\eta(J_\eta)}$ and 
$\hat{\xi}:\hat{V}\to\Delta\cap {\Bbb C}_{\xi(J_\xi)}$ respectively.
We will turn $\hol$ into a topological space by identifying it with
a subset of 
$\bran\times\bran$ by $\Ccal H\mapsto (\hat{U},0,\hat{\eta})\times(\hat{V},0,\hat{\xi}).$

As follows immediately from the definition,
\begin{prop}
\label{standard is epstein}
The holomorphic commuting pairs based on maps in the standard family
belong to $\hol$.
\end{prop}

Let us make a note of an important compactness property of $\Ccal E_s$
\begin{lem}
\label{compactness}
Let $s\in (0,1)$. The collection of normalized maps $g\in {\Ccal E}_s$
with $I=[0,1]$,  with marked domains $(U,0)$ 
is sequentially compact with respect to Carath\'eodory
convergence.
\end{lem}

\begin{pf}
Let $g_n:I=[0,1]\to J_n\subset \tl J_n$ be a normalized sequence in 
${\Ccal E}_s$. These maps extend to triple branched coverings
$g_n:G_n\to {\Bbb C}_{\tl J_n}$; and can be decomposed 
 (\ref{epstein-decomposition}) as $g_n=Q_n\circ h_n$ where 
$Q_n(z)=Q_{c_n}(z)=z^3+c_n$, and $h_n:I\to [0,b]$
is a univalent map $G_n\to \Delta(h_n)$.
By passing to a subsequence we can ensure that
$\tl J_n$ converge to an interval $\tl J$, and $c_n\to c$.
Then $(\Delta(h_n),0)$ Carath\'edory converges to a slit domain
$(\Delta,0)$.

Since $J_n$ is $s^{-1}$-commensurable with $I_n$, $h_n(1)$ is bounded.
As $g_n'(1)>s>0$, the derivatives $(h_n^{-1})'(h_n(1))$ are bounded
from above.
The points $h_n(1)$ stay away from the boundary of $\Delta$,
and it follows from Koebe Distortion theorem, that $\{h_n^{-1}\}$
form a normal family in $\Delta$.

Let us select a locally uniformly converging subsequence 
$h_n^{-1}$. Since $h_n(I)$ have bounded length,
the limit of this sequence is a non-constant, and therefore univalent,
function $h^{-1}:\Delta\to G$. Let $R_n:\Cbb D\to G_n$
be the normalized Riemann map, $R_n(0)=0$, $R'_n(0)>0$.
It can be decomposed as $R_n=h^{-1}_n\circ R'_n$,
where $R'$ is the normalized Riemann map $\Cbb D\to \Delta(h_n)$.
By the above, the maps $R_n$ converege locally uniformly
to the Riemann map $R:\Cbb D\to G$, which is equivalent to
Carath\'eodory convergence $(G_n,0)\to (G,0)$.

Finally, note that the convergence of $h_n^{-1}$ implies that
there is a point $c\in I$ where derivatives of $h_n$ are bounded
from above. It follows from Koebe theorem, that $\{h_n\}$
is a normal family in $G$, and hence $h_n\to h\equiv (h^{-1})^{-1}.$

\end{pf}

The above proof also yields:
\begin{lem}
\label{cubic-bounded-distortion}
For any $s\in (0,1)$, there exists a domain $O_s\supset [0,1]$,
such that for any $g\in {\Ccal E}_s$ with normalization $I=[0,1]$,
the univalent map $h$ in (\ref{epstein-decomposition})
is well-defined and has $K(s)$-bounded distortion in $O_s$.
\end{lem}

We will further refer to the above property by saying that a map
$g\in {\Ccal E}_s$ is {\it cubic up to bounded distortion}

The Epstein class of critical circle maps was first introduced by
Eckmann and Epstein \cite{EE} as an invariant subspace for the 
Renormalization operator, it was further shown in \cite{EE}
that this class contains a fixed point of $\cal R$ with golden
mean rotation number. It can be shown using real distortion
estimates, that renormalizations of any smooth circle map converge
to the Epstein class, and by recent work of de~Faria and de~Melo \cite{dFdM1},
the rate of convergence is geometric in $C^r$ metric. 
We shall only use a the following,  weaker statement:
\begin{lem}
\label{real-bounds}
Let $f\in C^r$, $(r\geq 3)$ be a critical circle map with an
 irrational rotation number.
Then the collection of real commuting pairs 
$(\wtl{{f}^{q_{m+1}}}|\wtl{{I}_{m}},\wtl{{f}^{q_{m}}}|\wtl{{I}_{m+1}})$
 is precompact in $C^r$ topology, and all its partial limits are 
contained in ${\cal E}_s$, for some universal constant $s>0$,
 independent on the original map $f$.
\end{lem}

In particular for a critical circle map $f\in\cal E$ 
there exists $\sigma>0$ such that all its renormalizations are contained 
in $\cal E_\sigma$. Moreover, the constant $\sigma$ can be chosen
independent on $f$, after  skipping first few renormalizations.
 
\begin{lem}
\label{unique fixed point}
Let $\zeta=(\eta,\xi)\in \EE$ be a critical commuting pair with 
$\rho(\zeta)=0$, which appears as
 a limit of a sequence $\{\zeta_n\}\subset \EE$
with $\rho(\zeta_n)\in \Cbb R\setminus \Cbb Q$.
 Then the map $\eta$ has a unique fixed point in
the interval $I_\eta$, which is necessarily parabolic, with
multiplier one.
\end{lem}
\begin{pf}
The existence of a fixed point in $I_\eta$ is clear. Since there are arbitrary
small real perturbations of $\eta$ without a fixed point on the real line,
any such point must be parabolic. 

Let us show that the fixed point is unique. 
We thank P. Jones for suggesting the following argument.
Assume that  $a$ and $b$ are two distinct fixed points of $\eta$ in $I_\eta$.
Since $\eta$ is in the
Epstein class, it has a well-defined inverse branch $\varphi$ in
${\Bbb C}_{[a,b]}$.
As $a$ and $b$ are fixed points, $\varphi([a,b])=[a,b]$.
Since $\varphi\ne \operatorname{Id}$, there exists $x_0\in [a,b]$ with $\varphi(x_0)\ne x_0$.
Without loss of generality, assume that $\varphi(x_0)>x_0$. Let $\dist_P(\cdot,\cdot)$ 
denote the Poincar\'e distance in ${\Bbb C}_{[a,b]}$. Set $\eps=\dist_P(\varphi(x_0),x_0)$.
By Schwarz Lemma, $\dist_P(\varphi(x),x)>\eps$ for all $x<x_0$.
On the other hand, $\dist_P(\varphi(x),x)<\operatorname{const}\cdot(x-a)^{-1}$.
Therefore, $\varphi'(a)=1+\eps'>1$, a contradiction.
\comm{
Elementary considerations imply that each of the points $a$, $b$ 
attracts some of the orbits of $\varphi$ in ${\Bbb C}_{[a,b]}$.
This contradicts Denjoy-Wolf Theorem.
\comm{
 Denote by $D$ the disk bounded
by $C$. Mapping $\varphi$  contracts the Poincar\'e metric on ${\Bbb C}_{[a,b]}$,
hence $D'\equiv\varphi(D)\subset D$. Viewing both domains $D$ and $D'$ as hyperbolic
manifolds, we see by symmetry that $[a,b]$ is a Poincar\'e geodesic joining $a$ to $b$
in both of these domains.
 Schwarz Lemma implies that $D'=D$, thus $\varphi$ is M\"obius.
Since $\varphi'(a)=\varphi'(b)=1$, $\varphi\equiv \operatorname{Id}$, which is impossible.}}
\end{pf}

A commuting pair $\zeta=(\eta,\xi)\in \EE$ will be called {\it parabolic} 
if the map $\eta$ has a unique fixed point in $I_\eta$, which has
a unit multiplier; this point will usually be denoted $p_\eta$.
Note, that
by virtue of its uniqueness, $p_\eta$ has to be globally attracting on one side
for the interval homeomorphism $\eta|_{I_\eta}$,
it is globally attracting on the other side under $\eta^{-1}$.

\medskip
\noindent
{\bf Complex bounds.}
The main analytic result of the paper \cite{Ya1} was establishing
 complex {\it a priori} bounds for Epstein 
critical circle maps in the following form:

\begin{thm}[\cite{Ya1}]
\label{complex bounds}
For any $s\in(0,1)$ there exists $N=N(s)$ and $D=D_r(0)$ such that
the following  holds: 

Let $f\in \EE_s$ be a critical circle map whose rotation number 
$\rho(f)$ has the continued fraction representation with at least
$N+1$ terms. Denote by $\Omega_m\ni 0$ the domain which triple covers
${\Bbb C}_{f^{q_{m+1}}(I_m)}$ under $f^{q_{m+1}}$. Then
\begin{equation}
\label{cubic estimate}
\frac{\dist(f^{q_{m+1}}(z),0)}{|f^{q_{m+1}}(I_m)|}+d>c\left(
\frac{\dist(z,0)}{|I_m|}\right)^3, \text{ for all }z\in \Omega_m,
\text{ with }f^{q_{m+1}}(z)\in D.
\end{equation}
 The constants $c$ and $d$ in the above inequality is universal.
\end{thm}

As an immediate consequence, for $m$ sufficiently large, we can
choose a Euclidean disc $\Delta\equiv D_{r_m}(0)$ around the origin, 
commensurable with $I_m$,
such that the preimages $U_m=f^{-q_{m+1}}(D_{r_m})\cap \Omega_m$
and $V_m=f^{-q_{m}}(D_{r_m})\cap \Omega_{m-1}$ are 
contained in a concentric disc with a smaller radius ${r'}_m$,
and moreover, $r_m/r'_m$ is greater than some fixed value $K>1$.
Thus the real commuting pair $(f^{q_{m+1}}|I_m,f^{q_m}|I_{m+1})$
extends to a holomorphic pair with range $\Delta$ and  definite modulus:
\begin{equation}
\label{extension}
{\cal H}:(f^{q_{m+1}}:U_m\to \Delta,f^{q_m}:V_m\to\Delta)
\end{equation}

For $\mu\in(0,1)$ let $\hol(\mu)$ denote the space of holomorphic commuting pairs
$\Ccal H:\Omega_{\Ccal h}\to \Delta_{\Ccal H}$, with $\mod (\Delta_{\Ccal H}\setminus\Omega_{\Ccal h})>\mu$,
$\min(|I_\eta|,|I_\xi|)>\mu$ and $\diam(\Delta_{\Ccal H})<1/\mu$.

\begin{lem}
\label{bounds compactness}
For each $\mu\in(0,1)$ the space $\hol(\mu)$ is sequentially pre-compact, with every limit point
contained in $\hol(\mu/2)$.
\end{lem}
\begin{pf}
Let $\{\Ccal H_n\}$ be a sequence of holomorphic pairs in $\hol(\mu)$, with ranges $\Delta_n$. 
Let $\hat{\eta}_n:\hat{U}_n\to \Delta_n$ and $\hat{\xi}_n:\hat{V}_n\to\Delta_n$ be the three
fold extensions of $\eta_n$, $\xi_n$.

By passing to a subsequence let us ensure that the Riemann maps $R_{(U_n,0)}$ converge locally
uniformly to a map $R$. It easily follows from Koebe Distortion Theorem that $R$ is non-constant and
thus a univalent map $R\equiv R_{(U,0)}$. 

The family $\hat{\eta}_n$ is normal in $U$ by Montel's Theorem, and passing to a subsequence again, 
we have $(\hat{U}_n,0,\hat{\eta}_n)\to (U,0,\hat{\eta})$. 
The convergence of $\hat{\xi}_n$ is ensured in the same fashion. Clearly, the resulting pair
has modulus greater than $\mu/2$, and thus is in $\hol(\mu/2)$.
\end{pf}

Let us define {\it germs of holomorphic commuting pairs} using the equivalence relation:
$\Ccal H\equiv\Ccal G$ if $K_{\Ccal H}$ and $K_{\Ccal G}$ coincide as sets, and $S_{\Ccal H}\equiv
S_{\Ccal G}$ on an open neighborhood of $K_{\Ccal H}$. The germ of a pair $\Ccal H$ will be denoted
$[\Ccal H]$.

\medskip

\noindent
{\bf Straightening theorems for holomorphic commuting pairs.}
The existence of  extensions to holomorphic pairs with
complex {\it a priori} bounds for high renormalizations
of Epstein circle maps, together with an analysis of the shapes of
the domains of these extensions was used in \cite{Ya1} to construct
quasiconformal conjugacies between high renormalizations of circle maps
with the same rotation number, with universally bounded dilatation.

Let $\zeta$ be at least $n$ times renormalizable
 critical commuting pair. For lack of a better term, let us say that
the pair of numbers $\tau_n(\zeta)=(r_{n-1},r_{n-2})$ forms the {\it history} of
the pair $\Ccal R^n\zeta$.
\begin{thm}[\cite{Ya1}]
\label{qc-conjugacy}
For each $s\in(0,1)$ there exists $N=N(s)$, such that the following holds.
Let $\zeta_1=(\eta_1,\xi_1)$ and $\zeta_2=(\eta_2,\xi_2)$ be two  critical 
commuting pairs in $\EE_s$ with irrational rotation numbers.
Assume that their $N$-th renormalizations have the same rotation number
and the same history.
 Then for every $n\geq N$ the renormalizations
${\cal R}^n\zeta_1$ and ${\cal R}^n\zeta_2$ extend to holomorphic commuting pairs
which are $K-$quasiconformally conjugate
with a universal dilatation bound $K$.
\end{thm}
For commuting pairs of bounded type this Theorem was proved by de~Faria \cite{dF1,dF2},
with ``$K$" depending on the bound on the rotation number.

One of the key points in quadratic-like renormalization theory is 
Douady-Hubbard Straightening Theorem, which claims that every quad\-ratic-like
map is conjugated to a quadratic polynomial via a quasiconformal homeomorphism
which is conformal on the Julia set (a so-called {\it hybrid equivalence}).
In renormalization theory of holomorphic commuting pairs,
 the role of quadratic polynomials is played by 
 holomorphic pairs  generated by standard maps 
(\lemref{standard-pairs}). Below we will establish certain rigidity properties
for these pairs, and prove a version of Straightening Theorem for holomorphic
commuting pairs.

\begin{lem}
\label{invariant differential}
Let $f_\theta$ be a map in the standard family. As before, let $A_\theta$
denote its lift to the plane, and $T(z)=z+1$.
Let $k\in\Cbb N\cup\{\infty\}$ denote the length of the continued fraction 
expansion of the rotation number $\rho(\theta)$.
Suppose that for some $n<k$, the critical commuting pair
\begin{equation}
\label{st-pair}
\zeta_{\theta,n}=(\eta=T^{-p_{n+1}}\circ A_\theta^{q_{n+1}},\; \xi=T^{-p_{n}}\circ A_\theta^{q_{n}})
\end{equation}
extends to a holomorphic commuting pair $\Ccal H:\Omega\to\Delta$,
where $p_n/q_n$ is the $n$-th convergent of $\rho(\theta)$.
Then any $\Ccal H$-invariant Beltrami differential $\mu$ with support in $K(\Ccal H)$
can be extended off $\Omega$ to a Beltrami differential $\hat\mu$ with support in 
$K(A_\theta)$ which is invariant
under $A_\theta$ and $T$ (that is $\hat\mu$ projects to a Beltrami differential on the
quotient cylinder $\Cbb C/\Cbb Z$, invariant under $f_\theta$).
\end{lem}

\begin{pf}
Set $g_n=T^{-p_{n}}\circ A_\theta^{q_{n}}:\Cbb C\to \Cbb C$.
Recall, that if $\rho(\theta)=[r_0,r_1,\ldots]$ then $q_{n+1}=r_n q_n+q_{n-1}$, with $q_0=1$,
$q_1=r_0$, and $p_{n+1}=r_n p_n+p_{n-1}$, with $p_0=0$, $p_1=1$.
Let $\mu$ be as above, and consider the Beltrami differential $\mu'$ obtained by
pulling $\mu$ back by various inverse branches of $g_{n+1}=T^{-p_{n+1}}\circ A_\theta^{q_{n+1}}$
and $g_n=T^{-p_{n}}\circ A_\theta^{q_{n}}$. Observe, that $\mu'$ is invariant under
$g_{n-1}=g_{n+1}\circ g_n^{-r_n}$. Arguing inductively, we see that $\mu'$ is
invariant under $g_0=A_\theta$ and $g_1=T\circ A_\theta^{r_0}$ and thus under $T$ as well.
\end{pf}

\begin{thm}[Straightening]
\label{straightening1}
Let $\zeta\in \EE$ be a critical commuting pair with an irrational rotation number.
 Then there exists $N$ such that for all $n>N$
the renormalization $\Ccal R^n(\zeta)$ is $K$-quasiconformally conjugate to a holomorphic
pair $\Ccal H$, whose underlying real commuting pair is $\zeta_{\theta,n}$ (\ref{st-pair})
for some $\theta\in[0,1)$, $n\in\Cbb N$.
The constant $K$ is universal, and the conjugacy is conformal on the filled Julia set $K(\Ccal H)$.
\end{thm}
\begin{pf}
The existence of the conjugacy is guaranteed by \thmref{qc-conjugacy}. Its conformality
on the filled Julia set follows from Lemmas \ref{invariant differential} and \ref{rigidity-irrational}.
\end{pf}

\begin{thm}[Straightening of limiting pairs]
\label{straightening2}
Let $\zeta\in \EE$ be a critical commuting pair with an irrational rotation number.
Assume that there is a sequence $\{n_k\}\subset\Cbb N$, such that some holomorphic
pair extensions $\Ccal H_k$ of $\Ccal R^{n_k}\zeta$ converge in $\hol$ to a pair 
$\Ccal H:\Omega_{\Ccal H}\to \Delta_{\Ccal H}$, which is parabolic. Assume in addition, that 
$\tau_{n_k}(\zeta)=(a,b)$.

Then we can find  $\theta=\theta(a,b)\in[0,1)$  such that 
the critical commuting pair $\zeta_{\theta,N}$ extends to a parabolic holomorphic 
commuting pair $\Ccal G:\Omega_{\Ccal G}\to\Delta_{\Ccal G}$
 with history $(a,b)$;
 there exists a $K$-quasiconformal map $\phi:\Delta_{\Ccal H}\to\Delta_{\Ccal G}$
which is a conjugacy: $\Ccal H=\phi^{-1}\circ\Ccal G\circ\phi$. The constant $K$ is 
universal, and $\phi$ is conformal on the filled Julia set $K(\Ccal H)$.
\end{thm}
\begin{pf}
The existence of $\Ccal G$, $K$-quasiconformally conjugate to $\Ccal H$ is guaranteed
by \thmref{qc-conjugacy} and compactness of quasiconformal maps.
Let $\psi:\Delta_{\Ccal H}\to\Delta_{\Ccal G}$ be a conjugacy.
Consider a new complex structure $\mu$ on $\Delta_{\Ccal G}$ given by the 
pullback $(\phi^{-1})^*\sigma_0$ on $K(\Ccal G)$ and the standard structure $\sigma_0$
elsewhere. 

By \lemref{invariant differential} the structure $\mu$ can be extended to a global
structure $\mu'$ invariant under $A_\theta$ and $T$. Let $ q:\Cbb C\to \Cbb C$ be the solution
of $\bar\partial {q}/\partial {q}=\mu'$ fixing $0,$ $1$, and $\infty$.
 By \lemref{rigidity of standard maps}, $ q\circ\Ccal G\circ q^{-1}=\Ccal G$.
The map $\phi=q\circ\psi$ is  the conjugacy with the desired properties.
\end{pf}
 
\comm{
\comm{
\medskip
\noindent
{\bf Some basic Teichm\"uller theory.} We will need some tools provided by Teichm\"uller theory, which we briefly recall below.
Let $X$ be a compact Riemann surface, and let $C$ be a closed at most countable
subset of $X$. The Teichm\"uller space of $X$ with boundary $C$ is defined to be the
set of isotopy classes of quasiconformal homeomorphisms of $X$ rel $C$. We denote it
by $T(X,C)$. The Teichm\"uller space is finite dimensional if $X$ has finite genus, and $C$ is a finite set. We will denote by $K_{\phi}(X)$ the maximal dilatation of 
a quasiconformal map $\phi$.
The {\it extremal map} in the isotopy class in $T(X,C)$
 is the one which has minimal quasiconformal dilatation.
A quasiconformal map is called a {\it Teichm\"uller map} if its beltrami differential is of the form $\mu=t\bar q/|q|$, where $q$ is a holomorphic quadratic differential
such that $||q||=\int_X|q|<\infty$, and $|t|<1$.

\begin{tet} Let $T(X,C)$ be a finite dimensional Teichm\"uller space. Then every isotopy class contains ann extremal map. Moreover, that map is unique and is a Teichm\"uller
map. If $T(X,C)$ is not finite dimensional, the extremal map exists but is not necessarily unique.
\end{tet}

Set $S=X\setminus C$ and let $K$ be any compact subset of $S$. Set 
$\hat K_\phi(S\setminus R)=\inf_{\psi\equiv\phi}K_\psi(S\setminus R)$. The {\it boundary dilatation} $H(\phi)$ is the direct limit of the numbers
$\hat K_\phi(S\setminus R) $ as $R$ increases to $S$.

\begin{sfmc} Let $\phi$ be a quasiconformal homeomorphism of $S$ to another surface,
and suppose $H(\phi)<\hat K_\phi(S\setminus R)$. Then the isotopy class of
$\phi$ rel $C$ contains a unique extremal map which is a Teichm\"uller map.
\end{sfmc}
}
}
\section{Parabolic  points, their perturbations, and parabolic renormalization}
\label{parabolic}

\noindent
{\bf General theory.}
      We begin with a brief review of the theory of parabolic bifurcations,
as applied in particular to an interval map in the Epstein class. For a more
comprehensive exposition the reader is referred to \cite{Do},
supporting technical details may be found in \cite{Sh}.
Fix a  map $\eta_0\in\Ccal E$ having a parabolic fixed point $p$
with unit multiplier.

\begin{thm}[Fatou Coordinates]
\label{Fatou-coord}
There exist topological discs $U^A$ and $U^R$, called {\it attracting}
and {\it repelling petals}, whose union is a 
punctured neighborhood of the parabolic periodic point $p$
such that
$$\eta_0(\bar U^A)\subset U^A\bigcup \{p\},\text{ and }
\bigcap_{k=0}^\infty \eta_0^{k}({ U}^A)=\{p\},$$

$$\eta_0(\bar U^R)\subset U^R\bigcup \{p\},\text{ and }
\bigcap_{k=0}^\infty \eta_0^{-k}({ U}^R)=\{p\},$$
where $\eta_0^{-1}$ is the univalent branch fixing $\zeta$.

Moreover, there exist injective analytic maps
$$\Phi^A:U^A\to{\Bbb C}\text{ and }\Phi^R:U^R\to \Bbb C,$$ 
unique up to post-composition by translations, such that
$$\Phi^A(\eta_0(z))=\Phi^A(z)+1\text{ and }\Phi^R(\eta_0(z))=\Phi^R(z)+1.$$
The Riemann surfaces $C^A=U^A/\eta_0$ and $C^R=U^R/\eta_0$ are
 conformally equivalent to the cylinder $\Bbb C/\Bbb Z$.
\end{thm}

We denote $\pi_A:U^A\to C^A$ and $\pi_R:U^R\to C^R$ the natural projections.
The quotients $C^A$ and $C^R$ are customarily referred to as {\it \'Ecalle-
Voronin cylinders}; 
 we will find it useful to regard these as Riemann
spheres with distinguished points $\pm$ filling in the punctures.
The real axis projects to {\it natural equators} $E^A\subset C^A$
and $E^R\subset C^R$.
Any conformal {\it transit homeomorphism} $\tau:C^A\to C^R$ fixing the ends $\pm$
is a translation in suitable coordinates.  Lifiting it produces a map
 $\bar \tau:U^A\to \Cbb C$ satisfying 
$$\tau\circ\pi_A=\pi_R\circ\bar\tau.$$
We will sometimes write $\tau\equiv \tau_\theta$, and $\bar\tau=\bar\tau_\theta$, where
$$\Phi^R\circ\bar\tau\circ(\Phi^A)^{-1}(z)\equiv z+\theta \mod(z).$$

If $z\in U^A\cap U^R$, we set $\Eps(\pi_R(z))=\pi_A(z)$. This {\it \'Ecalle-Voronin
transformation} extends to an analytic map of the neighborhoods $W(+)$, 
$W(-)$ of the two ends of $C_R$. Composing with a transit homeomorphism
yields an analytic dynamical system 
$F_\tau=\tau\circ{\Eps}:W(+)\cup W(-)\to C^R$ with fixed 
points $\pm$. The product of corresponding eigenvalues is clearly independent of $\tau$; noting that each of the components $W(\pm)$ is mapped onto the whole $C^R$ 
and applying Schwarz Lemma we conclude that 
\begin{equation}
\label{repelling ends}
|F'_\tau(+)|\cdot|F_\tau'(-)|>1.
\end{equation}

For an Epstein map $\eta$ in a sufficiently small neighborhood of $\eta_0$
the parabolic point splits into a complex conjugate pair of repelling
fixed points $p_\eta\in \Bbb H$ and $\bar p_\eta$ with multipliers
$e^{2\pi i\pm\alpha(\eta)}$. In this situation one may still speak of 
attracting and repelling petals:

\begin{lem}[Douady Coordinates]
\label{douady-coord}
There exists a neighborhood $U(\eta_0)\subset \Ccal E$ of the map $\eta_0$ 
such that for any $\eta\in U(\eta_0)$ with $|\arg \alpha(\eta)|<\pi/4$,
there exist topological discs $U_\eta^A$ and $U_\eta^R$ whose union is a 
neighborhood of $p$, and injective analytic maps
$$\Phi^A_\eta:U^A\to {\Bbb C}\text{ and }\Phi_\eta^R:U^R_f\to{\Bbb C}$$
unique up to post-composition by translations, such that
$$\Phi^A_\eta(\eta(z))=\Phi^A_\eta(z)+1
\text{ and }\Phi_\eta^R(\eta(z))=\Phi_\eta^R(z)+1.$$
The quotients $C^A_\eta=U^A_\eta/\eta$ and $C^R_\eta=U^R_\eta/\eta$
 are Riemann surfaces
conformally equivalent to ${\Bbb C}/{\Bbb Z}$.
\end{lem}

Note that the condition on arguments of
the eigenvalues of the repelling fixed points of $\eta$
is automatically satisfied for all maps in a sufficiently small neighborhood
of $\eta_0$ as follows from Holomorphic Index formula.

An arbitrary choice of real basepoints $a\in U^A$ and $r\in U^R$ enables us
to specify the Fatou and Douady coordinates uniquely, by requiring that
$\Phi^A(a)=\Phi^A_\eta(a)=0$, and $\Phi^R(r)=\Phi^R_\eta(r)=0$.
The following fundamental theorem first appeared in 
\cite{orsay-notes}:

\begin{thm}
\label{continuity-fatou}
With these normalizations we have
$$\Phi^A_\eta\to\Phi^A\text{ and }\Phi^R_\eta\to \Phi^R$$
uniformly on compact subsets of $U^A$ and $U^R$ respectively.

Moreover, select the smallest $n(\eta)\in\Cbb N$ for which 
$\eta^{n(\eta)}(a)\geq r$. Then 
$$\eta^{n(\eta)}(z)=(\Phi^R_\eta)^{-1}\circ T_{\theta(\eta)+K}\circ \Phi^A_\eta$$
wherever both sides are defined. In this formula $T_a(z)$ denotes the
translation $z\mapsto z+a$, $\theta(\eta)\in [0,1)$ is given by 
$$\displaystyle\theta(\eta)=1/\alpha(\eta)+
\underset{\alpha(\eta)\to\infty}{o(1)}\mod 1,$$ and the real
constant $K$ is determined by the choice of the basepoints $a$, $r$.
Thus for a sequence $\{\eta_k\}\subset U(\eta)$ converging to $\eta$, 
the iterates $\eta_k^{n({\eta_k})}$ converge locally uniformly if and only
if there is a convergence $\theta(\eta)\to \theta$, and the limit 
in this case is a certain lift of the transit homeomorphism $\tau_\theta$
for the parabolic map $\eta_0$.
\end{thm}

\comm{
A sequence $\{\eta_k\}$ converges to $\eta_0$
{\it tangentially} if $|\arg \_{f_k}|\to 0$. Tangential
convergence of a sequence is preserved under a quasiconformal change
of coordinates:

\begin{lem}
\label{tangential}
Consider two sequences of real-analytic maps:
$$f_n\to f_0,\text{ and }g_n\to g_0.$$
Assume that $f_n$ and $g_n$ are $K-$quasiconformally conjugate for
$n\geq 0$. Then if $f_n$ converge tangentially, so does $g_n$.
\end{lem}

\begin{pf}
Let ${\Bbb T}_{f_n}$ be the {\it quotient torus} of the periodic point
$\zeta_{f_n}$, defined as the quotient of an arbitrary linearizing neighborhood
of $\zeta_{f_n}$ under the action of $f_n^p$. One directly verifies
that 
$$\mod {\Bbb T}_{f_n}=\rho_{f_n}.$$
Appropriately define the quotient tori ${\Bbb T}_{g_n}$ with 
$\mod {\Bbb T}_{g_n}=\rho_{g_n}$. Observe that the quasiconformal conjugacies
induce $K$-quasiconformal homeomorphisms ${\Bbb T}_{f_n}\to {\Bbb T}_{g_n}$
preserving the preferred homotopy classes; the claim immediately
follows.
\end{pf}
}

\medskip
\noindent
{\bf Parabolic renormalization of  commuting pairs.}
Let $\zeta=(\eta,\xi)\in\EE$ be a parabolic commuting pair. 
Denote by $p\in I_\eta$ the parabolic fixed point of $\eta$. The pair $\zeta$
has zero rotation number and is, therefore, not renormalizable. In what follows
we will use the above discussed local theory of parabolic germs to describe
a {\it parabolic renormalization} construction for a parabolic pair $\zeta$,
which will naturally supplement the usual renormalization procedure.
A parallel construction for parabolic quadratic-like maps appeared in a 
paper of Douady and Devaney \cite{DD}.

Let $U^A$, $U^R$, $C^A$, and $C^R$ be as above. Fix an arbitrary transfer 
isomorphism $\tau:C^A\to C^R$ which preserves the equators 
(i.e. a rigid rotation). 
Let $N\geq 0$ be such that $\eta^N(\xi(0))\in C^A$. Fix a branch $\pi^{-1}_R$
and take the smallest $M$ for which 
$$\eta^M\circ\pi^{-1}_R\circ \tau\circ\pi_A\circ \eta^N\circ\xi(0)\in [0,\eta(0)].$$
Clearly, the composition 
$$\gamma\equiv \eta^M\circ\pi^{-1}_R\circ \tau\circ\pi_A\circ \eta^N\circ\xi$$
has a well-defined extension to the whole interval $[0,\eta(0)]$
which is independent of the choice of the branch of $\pi_R^{-1}$.

\begin{defn}
The { \it parabolic renormalization } of the commuting pair $\zeta=(\eta,\xi)$,
{\it  corresponding to the transit map $\tau$}, is the normalized pair
$$\pr_\tau\zeta=(\wtl{{\gamma}|{[0,\eta(0)]}},\wtl{{\eta}|{[\gamma(0),0]}}).$$
\end{defn}

\noindent
Again, we will refer to the non-rescaled commuting pair $(\gamma,\eta)$ as the 
 {\it pre-parabolic renormalization} $p{\pr}_\tau\zeta$ of the pair $\zeta$.

As an immediate consequence of continuity of Douady coordinates 
(\thmref{continuity-fatou}) we have the following:

\begin{lem}
\label{ren-converge}
Let $\zeta_k=(\eta_k,\xi_k)\in\EE$ 
be a converging sequence of renormalizable 
commuting pairs, $\zeta_k\to \zeta=(\eta,\xi)$,
with rotation numbers $\rho_k\to 0$. Assume that their
renormalizations also converge, ${\cal R}\zeta_k\to \tl\zeta$. Then
for some choice of transit map $\tau$ for the map $\eta$, we have
$${\cal P}_\tau\zeta=\tl\zeta.$$
\end{lem}

In the same way as for the  renormalization operator $\cal R$ (see  Remark \ref{injective1}) we have
\begin{rem}
\label{injective2}
The parabolic renormalization is injective. More precisely, for any $\zeta\in\crit$ there exists
at most one parabolic commuting pair $\zeta_{-1}\in\crit_\infty$ and a unique choice
of the transit map $\tau$ such that ${\cal P}_\tau\zeta_{-1}=\zeta$.
\end{rem}
\begin{pf}
Let $\zeta=(\eta,\xi)\in\crit_\infty$ be a parabolic pair, and let $\hat\zeta=(\gamma,\eta)=p{\cal P}_\tau\zeta$.
We can write $\gamma=\eta^M\circ\bar\tau\circ\eta^{N}\circ\xi$ for a choice of local lift $\bar\tau$ of the transit
map $\tau$ (note that all the maps in this composition commute). The number $N$ can be identified as follows.
Let $\bar\tau_1$ be a lift of an arbitrary transit map, then $N$ is the smallest natural number 
for which $\eta^{-(N+M)}\circ \bar\tau_1^{-1}\circ\gamma$ has a smooth inverse near $\eta(0)$.
The position of $\eta^M\circ\bar\tau\circ\eta^{N}\circ\xi(0)$ depends on $\tau$ monotonously,
which, in turn, specifies $\tau$.
\end{pf}

\section{Towers of holomorphic commuting pairs}
\begin{defn}
A {\it tower} of Epstein holomorphic commuting pairs is an element of 
a product space ${\hol}^N$, $N\leq \infty$,
$${\cal T}=({\cal H}_i)_{i=1}^{N},$$
for which the following holds. If $\zeta_i=(\eta_i,\xi_i)$
 denotes the real commuting pair
underlying $\Ccal H_i$, then either $\chi(\zeta_i)\ne \infty$, and 
$\tl\zeta_{i+1}=p{\Ccal R}\zeta_i$, or $\chi(\zeta_i)=\infty$, and 
$\tl\zeta_{i+1}=p{\pr}_{\tau_i}\zeta_i$, for some choice of the transit isomorphism $\tau_i$,
for all $i\leq N$.
\end{defn}
The pair ${\cal H}_1$ will be referred to as the {\it base pair} of the tower $\cal T$.
We shall denote by $\tw$ the space of all towers. For an element ${\Ccal T}\in \tw$
its {\it rotation number} $\rho({\Ccal T})$ is the sequence 
$$\rho({\Ccal T})=(\chi(\zeta_1),\chi(\zeta_2),\ldots)$$
The {\it domain} and {\it range} of the tower $\cal T$ is the domain $\Omega_1$ and range $\Delta_1$ of
the base pair ${\cal  H}_1$. 

The {\it dynamics} $F({\cal T})$ associated to the tower $\cal T$ is the collection of all finite compositions
$f=h_1\circ\cdots\circ h_k$, where $h_n$ is an element of ${\cal H}_i$ or a local lift
of some $\tau_i$. For a point $z\in \Omega_1$, the orbit of $z$ under $\cal T$
is the set ${\cal O}_{\cal T}(z)$ of images of $z$ under all elements of $F({\cal T})$ which are defined at $z$.
 We say that $z$ is an {\it escaping point } if 
${\cal O}_{\cal T}(z)\cap \Delta_1\setminus\Omega_1\ne \emptyset$. Non-escaping points form the 
{\it filled Julia set} $K({\cal T})$, its boundary is the {\it Julia set} $J({\cal T})$.

Two towers ${\cal T}_1=({\cal H}_i^1)_{i=1}^N$ and ${\cal T}_2=({\cal H}_i^2)_{i=1}^N$ with 
domains $\Omega_1$ and $\Omega_2$ 
are conjugate (quasiconformally, conformally, etc.) if there is a homeomorphism
$\phi$ defined on the neighborhood of the filled Julia set of ${\cal T}_1$
having the appropriate smoothness, satisying
$${\cal H}_i^2\circ\phi=\phi\circ{\cal H}_i^1\text{  for all  }i\leq N.$$
We will mostly be concerned with infinite towers (that is $N=\infty$).

Examples of towers are readily provided by holomorphic pairs with irrational
rotation numbers. More generally, let $\zeta_i=(\eta_i,\xi_i)$ 
be a sequence of Epstein commuting pairs with irrational rotation numbers, with
holomorphic extensions ${\cal H}_i\in\hol$. 
Suppose that the $n$-th pre-renormalizations $p{\cal R}^n\zeta_i$ have  holomorhic pair
extensions ${\cal H}_i^n$ which converge for every $n\leq N$ to some $\hat{\cal H}^i$.
 As follows from \lemref{ren-converge}, the sequence 
$\hat{\cal T}=(\hat{\Ccal H}_1,\hat{\Ccal H}_2,\ldots,\hat{\Ccal H}_,\ldots)$ forms a tower.
Such towers will be referred to as {\it limiting}.
\thmref{continuity-fatou} implies the following:
\begin{prop}
\label{geometric limit}
Let $\hat{\cal T}$ be a limiting tower as above. For any map $f\in F({\cal T})$, and an open subset 
$W$ compactly contained in the domain of definition of $f$, there exists a sequence $\{f_k\}$
of finite compositions of elements of ${\cal H}^n_i$ such that 
$f_k\rightrightarrows f$ in the domain $W$.
\end{prop}

We will  distinguish the case when $\zeta_i$ are pre-renormalizations of the same pair,
$\zeta_i=p{\cal R}^{n_i}\zeta$, by referring to $\hat {\cal T}$ as a {\it renormalization
limiting tower}.

Finally, we will say that a limiting tower $\hat{\cal T}$ is a {\it standard  tower} if its base pair
${\cal H}_1$ is a limit of pairs ${\cal G}_i$ with 
${\cal G}_i|_{\Bbb R}=\zeta_{\theta_i,n_i}$.

The main objective of this section is proving the following Combinatorial Rigidity
Theorem for renormalization limiting towers:
\begin{thm}
\label{combinatorial rigidity}
Any two renormalization limiting towers with the same rotation number are conjugate by
a quasiconformal homeomorphism, which is conformal on the filled Julia sets of the towers.
\end{thm}

The Theorem will follow from a rigidity statement for standard 
towers, to which it is reduced with the help of the Straightening Theorems.

\medskip
\noindent
{\bf Combinatorial rigidity of standard towers.} 
We begin the proof of \thmref{combinatorial rigidity} by establishing the following
fact about standard  towers:
\begin{thm}
\label{one tower}
Let ${\cal T}=({\cal H}_i)$ and $\tl{\cal T}=(\tl{\cal H}_i)$ be two standard 
 towers. Assume that the base real commuting pairs 
$\zeta_1={\cal H}_1|_{\Bbb R}$ and $\tl\zeta_1=\tl{\cal H}_1|_{\Bbb R}$ and
the rotation numbers $\rho({\cal T})$ and $\rho(\tl{\cal T})$ coincide.

Then ${\cal H}_i|_{\Bbb R}=\tl{\cal H}_i|_{\Bbb R}$ for all $i\in\Bbb N$.
\end{thm}

The claim of the \thmref{one tower} will evidently follow from 
\begin{lem}
\label{one transit}
Let $\zeta\in\cal E$ be a parabolic commuting pair.
Let $\{r_i\}_1^m$, $m\leq\infty$ be a sequence of positive integers,
ending in $\infty$ if $m$ is finite, and
set $\rho=[r_1,\ldots]$.

\noindent
There exists a unique choice of the transit map $\tau$
for which 
$\rho({\pr}_\tau\zeta)=[r_1,\ldots]$, and in the case when $\rho\in\Bbb Q$ (that is $m<\infty$)
${\cal R}^{m-1}({\pr}_\tau\zeta)$ is a parabolic pair.
\end{lem}
\begin{pf}
We need to introduce some further notation. Let $p$ denote the parabolic fixed point
of $\eta$, and let $C^A$, and $C^R$ be as in 
\secref{parabolic}. The dynamics of the commuting pair $(\eta,\xi)$ induces a 
natural return map of the equators of the \'Ecalle cylinders $r:E^R\to E^A$.
For a choice of the transit isomorphism $\tau:C^A\to C^R$ fixing the equators
the composition 
$$f_\tau=\tau\circ r:E^A\to E^A$$
is a critical circle map. 
It is easily seen that $f_\tau$ is topologically conjugate to the circle map
associated to ${\pr}_\tau\zeta$ (see \propref{rotation number}), thus
$$\rho(f_\tau)=\rho({\pr}_\tau\zeta).$$
Note that the dependence $\tau\mapsto f_\tau(x)$ is a monotone map ${\bf T}\to\bf T$ for any $x\in T$.
The standard considerations then imply that $\psi:\tau\mapsto\rho(f_\tau)$ is a non-decreasing
degree one continuous map of the circle; moreover $\psi^{-1}(\rho)$ is a single point for each 
irrational $\rho$.

It remains therefore to consider the case when $\rho=p/q$ is rational. The existence of a $\tau$ satisfying
the conditions of the Lemma follows from \lemref{ren-converge}, let us establish its uniqueness.
Let $F_\tau:{\Bbb R}\to {\Bbb R}$ be a lift of $f_\tau$ having singularities
at integer points, 
with $F(0)\in [0,1)$. The map $f_\tau$
has rotation number $p/q$ if and only if $F_\tau^q(x)=x+p$ for some $x\in\Bbb R$.
If $\tau$ satisfies the conditions of the Lemma, the pair ${\pr}_\tau\zeta$ has a unique periodic 
orbit. The uniqueness of the orbit implies that  $F(x)\leq x$ for all $x$, or
$F(x)\geq x$ for all $x$. Monotone dependence of $f_\tau$ on $\tau$ now implies
that the value of $\tau$ realizing each of these two possibilities is unique.
Moreover, since the family $f_\tau$ contains no rigid rotations, the two
values are distinct.
Finally we note that by our conventions the continued fraction expansions of $\rho({\pr}_\tau\zeta)$
in these two situations will differ.

\end{pf}

Theorems \ref{rigidity of standard maps} and \ref{one tower} have the following important 
consequence:
\begin{thm}
\label{no deformations}
A standard tower admits no non-trivial quasiconformal deformations
entirely supported on the filled Julia set of the base pair.
\end{thm}

\medskip
\noindent
{\bf Proof of \thmref{combinatorial rigidity}.}
Let ${\cal T}=({\cal H}_1,{\cal H}_2,\ldots)$ be a limiting tower for a sequence
${\cal T}_i$ of towers based on the extensions ${\cal Z}_i$ of the 
renormalizations ${\cal R}^{k_i}\zeta$ of 
a commuting pair $\zeta\in \EE$. By \thmref{straightening1}, there exists a 
sequence of standard pairs ${\cal G}_i$  quasiconformally
conjugate to ${\cal Z}_i$ via a homeomorphism $\phi_i$ 
with a uniform dilatation bound. We can select 
the sequence ${\cal G}_i$ in such a way that their underlying pairs 
$\zeta_{\theta_i,n_i}$ converge to a pair $\hat\zeta$. There is no canonical choice
for $\hat\zeta$, but we may associate a fixed $\hat\zeta$ to all renormalization
limiting towers ${\cal T}$ with the same $\rho({\cal T})$.
Using \thmref{complex bounds}, \lemref{bounds compactness} and
compactness properties of quasiconformal maps, we may select a subsequence
$\phi_{k_i}$ converging to a quasiconformal homeomorphism $\phi$,
conjugating ${\cal H}_1:\Omega_1\to\Delta_1$ to a 
pair ${\cal G}$ with ${\cal G}|_{\Bbb R}=\hat\zeta$:
$${\cal H}_1\circ\phi=\phi\circ{\cal G}.$$
Let us denote by $\cal S$ the standard  tower with base $\cal G$,
it is unique by \thmref{one tower}.

We define a new Beltrami differential $\mu$ on $\Delta_1$ by first setting $\mu=\phi^*\sigma_0$ for $z\in\Delta_1\setminus\Omega_1$.
For a point $z$ in $\Omega_1$ whose orbit enters $\Delta_1\setminus\Omega_1$
we set $\mu=g^*\phi^*\sigma_0$, where $g$ is a composition of the elements of $\cal T$
mapping $z$ to $\Delta_1\setminus\Omega_1$. On the filled Julia set $K({\cal T})$ we
set $\mu=\sigma_0$.
Defined in this way, $\mu$ is invariant under the tower $\cal T$. Let 
$h:\Delta_1\to\Bbb C$ be the 
quasiconformal homeomorphism integrating $\mu$, $\mu=h^*\sigma_0$. Then
$S'\equiv h\circ{\cal T}\circ h^{-1}$ is a tower. The base pair of $S'$
is a quasiconformal deformation $h\circ\phi^{-1}\circ {\cal G}\circ\phi\circ h^{-1}$
of ${\cal G}$ entirely supported on its filled Julia set. It is thus equal to 
${\cal G}$ on a neighborhood of the filled Julia set. The rotation number of the
standard tower $S'$ is equal to $\rho({\cal T})$ and therefore, by 
\thmref{combinatorial rigidity}, $S'=S$.
Q.e.d.

\subsection{Bi-infinite towers}
A bi-infinite tower is an element of the product space $\hol^{\Bbb Z}$
$${\cal T}=(\ldots,{\cal H}_{-n},\ldots {\cal H}_0,\ldots,{\cal H}_n,\ldots)$$
with the following properties. For each $n\in\Bbb Z$ setting $\zeta_i=(\eta_i,\xi_i)={\cal H}_i\cap\Bbb R$
we either have $\chi(\zeta_i)\ne \infty$ and $\tl\zeta_{i+1}=p{\cal R}\zeta_i$;
or $\chi(\zeta_i)=\infty$ and $\tl \zeta_{i+1}=p{\pr}_{\tau_i}\zeta_i$ for some choice of the transit map
$\tau_i$.
The {\it rotation number} of $\cal T$ is naturally defined to be the bi-infinite sequence 
$$\rho({\cal T})=(\ldots,\chi(\zeta_{-i}),\ldots,\chi(\zeta_0),\ldots,\chi(\zeta_i),\ldots).$$

Examples of bi-infinte towers can be constructed from renormalizations of Epstein commuting pairs
with irrational rotation numbers as follows. Let $\zeta_k\in\EE$ with $\rho(\zeta_k)\in{\Bbb R}\setminus\Bbb Q$.
Let 
$${\cal T}_k=({\cal H}_{-i_k}^k,{\cal H}_{-i_k+1}^k,\ldots,{\cal H}_0^k,{\cal H}_1^k,\ldots)$$
be a forward tower of holomorphic pair extensions of the renormalizations of $\zeta_k$
$${\cal H}_{j-i_k}^k\cap {\Bbb R}=\lambda_k^{-1}\circ{\cal R}^j\zeta_k\circ\lambda_k$$
where the homothety $\lambda_k$ is chosen so that $\tl{\cal H}^k_0={\cal H}^k_0$.
If $\zeta_k$ are selected so that $i_k\to\infty$ and also 
${\cal H}_n^k\underset{k\to\infty}\longrightarrow{\cal H}^k$ for all $n\in\Bbb Z$ then
the sequence ${\cal T}=({\cal H}^k)^{\infty}_{-\infty}$ is a bi-infinite tower.
If in addition $\zeta_k={\cal R}^{n_k}\zeta$ for some $\zeta\in\EE$ then ${\cal T}$ is
a {\it limiting renormalization tower}.

The main result of this paper transforms into the following uniqueness theorem for bi-infinite
limiting renormalization towers:
\begin{thm}
\label{uniqueness}
Let ${\cal T}^1=({\cal H}^1_k)_{-\infty}^{\infty}$ and ${\cal T}^2=({\cal H}^2_k)_{-\infty}^{\infty}$
be two limiting renormalization towers with the same rotation number. Then
$${\cal H}^1_k\cap{\Bbb R}=\lambda^{-1}\circ{\cal H}^2_k\circ\lambda$$
for some real homothety $\lambda$ for all $k\in\Bbb Z$.
\end{thm}

The proof of \thmref{uniqueness} is broken into two steps. We first establish
\begin{thm}
\label{cbit}
Let ${\cal T}^1=({\cal H}^1_k)_{-\infty}^{\infty}$ and ${\cal T}^2=({\cal H}^2_k)_{-\infty}^{\infty}$
be two limiting renormalization towers with the same rotation number. Then
there exists a quasiconformal homeomorphism $\phi:{\Bbb C}\to{\Bbb C}$
conjugating the towers:
$${\cal H}^1_k=\phi^{-1}\circ{\cal H}^2_k\circ \phi,\text{ for all }k\in\Bbb Z.$$
\end{thm}
\begin{pf}
By \thmref{combinatorial rigidity} the towers 
$${\cal T}^1_n=({\cal H}^1_n,{\cal H}^1_{n+1},\ldots)\text{ and }
 {\cal T}^2_n=({\cal H}^2_n,{\cal H}^2_{n+1},\ldots)$$
are quasiconformally conjugate by a homeomorphism $\phi_n:\Delta_n^1\to\Delta_n^2$
of the domains of ${\cal H}^1_n$ and ${\cal H}^2_n$ respectively,
whose dilatation is bounded by a universal constant.
By complex bounds, the domains $\Delta_n^i$ fill out the plane;
choosing a diagonal subsequence of maps $\phi_{n_i}$ converging in each
$\Delta_n$ we obtain a limiting quasiconformal mapping $\phi:{\Bbb C}\to {\Bbb C}$ with the 
required properties.
\end{pf}

We now proceed to formulate
\begin{thm}
\label{bdl}
A limiting renormalization bi-infinite tower admits no nontrivial invariant Beltrami differentials.
\end{thm}
\thmref{bdl} readily implies \thmref{uniqueness}, since the Beltrami differential of the
quasiconformal conjugacy $\phi$ produced by \thmref{cbit} is invariant under the tower
${\cal T}^1$, and so equal to zero almost everywhere. Thus the conjugacy $\phi:{\Bbb C}\to{\Bbb C}$
is conformal and therefore a homothety.
The proof of \thmref{bdl} is based on hyperbolic metric expansion techniques developed
by McMullen in \cite{McM2} and used by Hinkle in the context parallel to ours in \cite{Hinkle}.
It will occupy most of the remainder of the paper.

\subsection{Expansion of the hyperbolic metric}
We fix a bi-infinite renormalization limiting tower ${\cal T}=({\cal H}_k)_{-\infty}^{\infty}$,
${\cal H}_k\cap{\Bbb R}=\zeta_k=(\eta_k,\xi_k)$.
We will refer to  $h\in F({\cal T})$ as a {\it map of  level $k$ } if $h$ is either an element of
the holomorphic commuting pair ${\cal H}_k$, that is one of the maps $\eta_k$, $\xi_k$, or $\eta_k\circ\xi_k$,
or a lift of a  transit homeomorphism $\tau$ associated to the parabolic point of $\eta_k$.
The domain of $h$, denoted by $D(h)$, will in the first case
denote the domain of the extension of $h$ to a degree three proper map
onto the plane; 
and in the case when $h$ is a lift of a transit map
it will be the domain of the maximal extension of $h$ provided by the Fatou coordinates.

Given a hyperbolic Riemann surface $X$ we shall denote by $d_X(\cdot,\cdot)$ the Poincar\'e 
distance on $X$. For a differentiable map $f:X\to Y$ of hyperbolic Riemmann surfaces,
$||f'(x)||_{X,Y}$ will stand for the norm of the derivative with respect to the
hyperbolic metrics; we will simply write $||f'(x)||$ if $X=Y={\Bbb C}\setminus{\Bbb R}$.
The expansion properties of bi-infinite towers rely on the following lemma formulated by McMullen
in \cite{McM2}:
\begin{lem}
\label{inclusion expansion}
There exists a continuous and increasing function $C(s)<1$ with $C(s)\to 0$ as $s\to 0$
such that for  the inclusion $\iota$ of a hyperbolic Riemann surface $X$ into a hyperbolic 
Riemann surface $Y$,
$$||\iota'(z)||_{X,Y}<C(s),$$
where $s=d_Y(x,Y\setminus X)$.
\end{lem}

The following estimate for the variation of the norm of the derivative of an analytic map
of hyperbolic Riemann surfaces is a consequence of Koebe Distortion Theorem:
\begin{lem}[Corollary 2.27\cite{McM1}]
\label{variation distortion}
Let $f:X\to Y$ be an analytic map between hyperbolic Riemann surfaces with nowhere vanishing derivative.
Then for $z_1, z_2\in X$ we have
$$||f'(z_1)||_{X,Y}^{1/\alpha}\leq||f'(z_2)||_{X,Y}\leq||f'(z_1)||_{X,Y}^{\alpha},$$
where $\alpha=exp(Kd_X(z_1,z_2))$ with a universal constant $K>0$.
\end{lem}

We now apply the above expansion principles to the setting of limiting renormalization towers:
\begin{lem}
\label{expansion towers}
Let ${\cal T}=({\cal H}_1,{\cal H}_2,\ldots)$ be a limiting renormalization tower, $h\in F({\cal T})$.
As before, denote by $\Omega_n=U_n\cup V_n\cup D_n$ and $\Delta_n$ the domain and range of the 
holomorphic pair ${\cal H}_n$.
We have the following:
\begin{itemize}
\item[(I)] $||h'(z)||>1$ for any $x\in D(h)$;
\item[(II)] There exists a universal constant $C>1$ such that $||h'(z)||>C$ if $z\in\Omega_n$ and
$h$ is an  element of level $n$ for which $h(x)\in\Delta_n\setminus\Omega_n$.
\end{itemize}
\end{lem}

\begin{pf}
Let us consider a holomorphic commuting pair $\cal H$ 
which is an extension of a high renormalization
${\cal R}^n\zeta$ for $\zeta\in\EE$ with $\rho(\zeta)\in{\Bbb R}\setminus{\Bbb Q}$.
Following the earlier convention, we write
$${\cal H}=(\eta=f^{q_{n+1}}|_U,\xi=f^{q_n}|_V).$$
The pair $\cal H$ should be viewed as a small perturbation of an element of $\cal T$.
The claim of the Lemma will follow by continuity, once we establish the properties (I) and (II)
for the elements of $\cal H$.

Consider first the map $\eta:U\to\Delta\cap {\Bbb C}_{\eta(J_U)}$. By the convention we have made,
$D(\eta)$ will denote the domain of its three-fold extension $\eta:D(\eta)\to {\Bbb C}_{\eta(J_U)}$.
Set $D\equiv D(\eta)\setminus\eta^{-1}{\Bbb R}\subset{\Bbb C}\setminus{\Bbb R}$, and
$\iota:D\hookrightarrow {\Bbb C}\setminus{\Bbb R}$. By \lemref{inclusion expansion} 
$$||\iota'(z)||_{D,{\Bbb C}\setminus{\Bbb R}}<C(s)<1,\text{ where }s=\dist_D(z,\partial D).$$
On the other hand, $\eta:D\to{\Bbb C}\setminus{\Bbb R}$ is a local isometry with respect to
the hyperbolic metrics, that is $||\eta'(z)||_{D,{\Bbb C}\setminus{\Bbb R}}=1$.
Thus, $||\eta^{-1}||>C(s)^{-1}>1$, which proves (I) for the map $\eta$.

Let us now establish (II), that is show that for $z\in U$ with $\eta(z)\in\Delta\setminus\Omega$,
$||\eta'(z)||>C>1$ for some fixed $C$. 
Note that $U\cap {\Bbb R}=J_U=[0,f^{q_n-q_{n+1}}(0)]$, and $\eta([0,f^{q_n-q_{n+1}}])=[q_{n+1},q_n]$,
which is well inside $\Omega$. This, real {\it a priori} bounds, and Koebe Theorem, imply 
that for some universal $\delta>0$ 
$$f^{q_{n+1}}(U_\delta([0,f^{q_n-q_{n+1}}(0)]))\subset \Omega.$$
Let us now recall the analysis of the shape of the domain $U$ carried out in \cite{Ya1}.
We summarize the relevant results of \cite{Ya1} as follows. 
There exists a universal constant $\eps>0$, such that for any point $z\in U\setminus U_\delta([0,f^{q_n-q_{n+1}}(0)])$ either $\operatorname{Im}z>\eps|[0,f^{q_n-q_{n+1}}(0)]|$
or $z$ is contained in a possibly empty set $W'$ with the following properties.
The set $W'$ is non-empty only if the iterate
 $f^{q_n}|_{[f^{q_{n-1}}(0),0]}$ is a sufficiently small perturbation of a parabolic map.
In this case $W'$  is an $\Bbb R$-symmetric topological disk whose inner diameter is commensurable
with $[0,f^{q_n-q_{n+1}}(0)]$; and
 there exists a fundamental
domain $C$ for the Douady coordinate of $f^{q_n}$ such that the univalent image $W=f^{q_{n-1}}(W')$
is obtained by repeated translation $C_0=C$, $C_1$, $C_2,\ldots,C_k=C'$ of $C$ under $f^{q_n}$.
The interiors of all the (crescent-shaped) regions $C_i$, except for the last one, are disjoint from $C$.
The intersection $(f^{q_{n-1}}(\partial U)\cap W)\cap \Bbb H$ is a connected simple curve 
$\Gamma$ (a ``horocycle")
obtained by translation of a fundamental arc $\gamma\subset C$ under $f^{q_n}$.
Note that by compactness of the Epstein class the equatorial annulus 
which the projections of $\Gamma$ and $-\Gamma$ cut out on the Douady cylinder of $f^{q_n}$ has
definite modulus.

Now let $z$ be a point of $U$ with $\eta(z)=f^{q_{n+1}}(z)\notin \Omega$.
If $z\notin W'$, then $\dist_{{\Bbb C}\setminus {\Bbb R}}(z,\partial D)<s_0$ with some universal
bound $s_0$ and we have (II) for $z$.
Otherwise, $w=f^{q_{n-1}}(z)$ is contained in one of $C_i$. Let $m$ denote the iterate
$c'=f^{mq_n}(c)\in C'$. Compactness considerations on the shape of $D(f^{q_n})$ now imply
$$\dist_{{\Bbb C}\setminus{\Bbb R}}(c',\partial D(f^{q_n}))>s_1>0$$
for some universal $s_1$. Hence $||Df^{q_n}(c')||>C>1$. Combining this and property (I) with the 
chain rule for decomposition
$f^{q_{n+1}}=f^{q_n}\circ\cdots\circ{f^{q_n}}\circ f^{q_{n-1}}$ we have $||\eta'(z)||>C>1$.

The statement (I) is proved in the same way for $\xi$ and $\eta\circ \xi$; (II) is obvious for 
$\eta\circ\xi$ and is proved in an identical fashion for $\xi$. 
The statements transfer to limiting renormalization towers by the standard considerations of continuity.

\end{pf}

\subsection{The structure of the filled Julia set of a limiting renormalization tower.}
In this section we prove the following theorem:
\begin{thm}
\label{structure}
The filled Julia set $K({\cal T})$ of a limiting renormalization tower $\cal T$ has no interior,
$K({\cal T})=J({\cal T})$. Moreover, repelling periodic orbits of maps in $F({\cal T})$ are
dense in $J({\cal T})$.
\end{thm}

\begin{lem}
\label{only parabolics}
let ${\cal T}$ be a limiting renormalization tower. Then every non-repelling periodic
orbit of ${\cal T}$ contains the parabolic point of some ${\zeta}_n$ in ${\cal T}$.
\end{lem}
\begin{pf}
Since the claim of the lemma is certainly true for periodic orbits in the real line,
let us assume that $z_1,\ldots,z_n$ is a periodic orbit of a map $h\in F({\cal T})$
disjoint from ${\Bbb R}$. By \lemref{expansion towers}, 
$$||Dh^{\circ n}(z_1)||>1,$$
which implies that the orbit $z_1,\ldots,z_n$ is repelling.
\end{pf}

Let us first establish the density of the repelling orbits in $J({\cal T})$,
\begin{lem}
\label{periodic dense}
Repelling periodic orbits of maps in $F({\cal T})$ are dense in $J({\cal T})$. Moreover, for any point $z\in
J({\cal T})$ there exists an element $h\in F({\cal T})$ mapping a neighborhood of $z$ onto the domain 
of $\cal T$.
\end{lem}
\begin{pf}
Let $W$ be an open neighborhood of a point $p\in J({\cal T})$. Clearly, there is an
element $h\in F({\cal T})$ such that $h(W)$ intersects the Julia set of a holomorphic pair
${\cal H}\in{\cal T}$. Let us perturb ${\cal H}$, if necessary, to a pair ${\cal H}'$
with an irrational rotation number, so that $W$ still intersects $J({\cal H}')$.
By reducing $W$, if necessary, we may assume that $W\cap {\Bbb R}=\emptyset$.
Denote by ${\cal G}$ a holomorphic extension of the standard pair $\zeta_{\theta,0}$
quasiconformally conjugate to ${\cal H}'$. The conjugacy corresponds to $W$ a 
neighborhood $\tl W$ intersecting the Julia set of the Arnold map ${A}_\theta$. 
Thus there is an iterate of $A_\theta$ mapping $\tl W$ over the whole domain of $\cal G$.
The considerations of convergence imply the existence of an element $h_W\in F({\cal T})$
mapping $W$ onto the domain of $\cal T$.
The existense of a repelling periodic point of $\cal T$ in $W$ follows from \lemref{only parabolics}.
\end{pf}

To prove \thmref{structure} it remains to show that $K({\cal T})$ has no interior.

Let us call a component $U$ of the interior of $K({\cal T})$ {\it wandering} if it is
disjoint from its forward images under maps of $F({\cal T})$. A non-wandering
component will be called {\it periodic}. 
A modification of the argument used in the above lemma implies that a component of $\overset{\circ}{K}({\cal T})$ 
does not intersect the Julia set of any of the holomorphic pairs $\cal H$ forming $\cal T$.
This implies that every map $h\in F({\cal T})$ defined on a subset of $U$ is  defined on 
all of $U$. Thus for a periodic component $U$ there exists an $h\in F({\cal T})$ with $h(U)=U$.
The disjointness from the Julia sets also implies that $U\cap {\Bbb R}=\emptyset$
for every $U\subset \overset{\circ}{K}({\cal T})$.

\begin{lem}
\label{no wandering}
The filled Julia set of a limiting renormalization tower has no wandering interior components.
\end{lem}
\begin{pf}
Assume the contrary. Note that since interior components of $K({\cal T})$ do not intersect the real
line, $h|_U$ is univalent for any  $U\subset \overset{\circ}{K}({\cal T})$ and $h\in F({\cal T})$ defined
on $U$. The appropriate case of Sullivan's no wandering domains Theorem can be directly
translated to our setting to construct non-trivial quasiconformal deformations of the tower $\cal T$
entirely supported on the grand orbit of the wandering domain. This contradicts the statement
of \thmref{no deformations}.
\end{pf}

We recall the following fundamental principles of dynamics on hyperbolic Riemann surfaces (see e.g. \cite{Lyubich-survey}).
\begin{prop}
\label{possibilities}
let $h:U\to U$ be an analytic transform of a hyperbolic Riemann surface. Then
one of the following four possibilities holds:
\begin{itemize}
\item[(I)]$h$ has an attracting or superattracting fixed point in $U$ to which all points converge;
\item[(II)] all orbits tend to infinity;
\item[(III)] $h$ is conformally conjugate to an irrational rotation of a disk, a punctured
disk or an annulus;
\item[(IV)] $h$ is a conformal homeomorphism of finite order.
\end{itemize}
\end{prop}
The next proposition expands on the second possibility:
\begin{prop}\label{possibility2}
Let $h$ be an analytic transform of a hyperbolic domain $U\subset \hat{\Bbb C}$,
continuous up to the boundary of $U$. Suppose that the set of fixed point of $h$ on $\partial U$ is totally disconnected. If the second possibility of \propref{possibilities} occurs, there is a fixed point $\alpha\in\partial U$
such that $h(z)\to\alpha$ for every $z$.
\end{prop}

The following lemma rules out the existense of periodic components for a limiting
renormalization tower, thus completing the proof of \thmref{structure}.
\begin{lem}
\label{no periodic}
The filled Julia set of a limiting renormalization tower has no periodic components.
\end{lem}
\begin{pf}
Assume the contrary. Let $U\subset \overset{\circ}{K}({\cal T})$ be a periodic 
component, and denote by $h$ the element of $F({\cal T})$ fixing $U$.
Let us consider the possibilities of \propref{possibilities}. The expansion
properties of $h$ (\lemref{expansion towers}) rule out the cases
(III) and (IV). The case (I) is ruled out by \lemref{only parabolics}.
The only remaining possibility is (II). Since $h$ is not the identity map,
the fixed points of $h$ in $\partial U$ are isolated, and 
\propref{possibility2} implies that all orbits in $U$ converge to a fixed point $p$
in $\partial U$. Again applying \lemref{only parabolics} we see that $p\in \Bbb R$
is the parabolic fixed point of one of the commuting pairs forming $\cal T$.

Denote by $\tl U$ the projection of the domain $U$ onto the repelling Fatou cylinder
$C^R$ of $p$. The domain $\tl U$ contains the ends $\pm$ of the cylinder. On the
other hand, the return map of $C^R$ is repelling at the
ends (\ref{repelling ends}). This implies,  that
$\tl U$ contains a preimage of the equator ${\Bbb R}/{\Bbb Z}$, and thus 
$U$ intersects $F({\cal T})^{-1}({\Bbb R})$. The last statement contradicts the assumption
that $U$ is an interior component of $K({\cal T})$.
\end{pf}

\subsection{Quasiconformal rigidity of limiting towers: proof of \thmref{bdl}}
\begin{thm}
\label{rigidity forward}
A limiting renormalization tower supports no invariant Beltrami differentials on its filled Julia 
set.
\end{thm}

Before giving the proof of the theorem let us state the following
\begin{lem}[cf. Lemma 1.8 \cite{Lyubich-survey}]
\label{disconnected}
Let $\cal T$ be a limiting renormalization tower. The group $G$ of homeomorphisms of $K({\cal T})$
which commute with all maps $h\in F({\cal T})$ is totally disconnected.
\end{lem}
\begin{pf}
Let $\phi\in G$ be a map in the connected component of the identity. Suppose $z_0\in K({\cal T})$ is
a repelling periodic point with $h(z_0)=z_0$ for some $h\in F({\cal T})$. Since the solutions
of $h(z)=z$ are isolated, $\phi$ fixes $z_0$. The claim now follows from  density of repelling cycles,
\thmref{structure}.
\end{pf}

\medskip
\noindent
{\it Proof of \thmref{rigidity forward}.} In view of \thmref{combinatorial rigidity} it is enough to
prove the statement for standard limiting towers. Let $\cal T$ be such tower, and let $\mu$ be
an invariant Beltrami differential with support in $K({\cal T})$. Let $\sigma_t$ denote
the  complex structure in the plane given by the standard structure $\sigma_0$ on 
${\Bbb C}\setminus K({\cal T})$, and by $t\mu$ on $K({\cal T})$, where $0\leq t\leq 1$.
Let $\phi_t$ be the normalized solution of Beltrami equation $\phi_t^*\sigma_0=\sigma_t$.
By \thmref{no deformations} the towers $\phi_t\circ {\cal T}\circ \phi_t^{-1}$ coincide with $\cal T$.
Thus each map in the continuous family $\phi_t$ commutes with all $h\in F({\cal T})$. As $\phi_0\equiv \operatorname{Id}, $
\lemref{disconnected} implies that $\phi_t|_{K({\cal T})}\equiv \operatorname{Id}$ for all $t$.
On the other hand, $\phi_t$ is conformal off $K({\cal T})$. Bers' sewing lemma implies that $\phi_t$
is conformal in the whole plane, and thus $\phi_t\equiv\operatorname{Id}$.
This implies that $\mu$ vanishes almost everywhere on $K({\cal T})$. $\;\Box$

\medskip 
For a bi-infinite  tower ${\cal T}=({\cal H}_k)^{\infty}_{-\infty}$ let ${\cal T}^n$
denote the truncated tower $({\cal H}_k)^{\infty}_{-n}$. Also denote by $\Omega_n$ and
$\Delta_n$ the domain and range of the pair ${\cal H}_n$. 
\begin{thm}
\label{fill the plane}
Let ${\cal T}=({\cal H}_k)^{\infty}_{-\infty}$ be a limiting renormalization
bi-infinite tower. Then
$$\displaystyle\lim_{n\to-\infty}J({\cal T}_n)=\Bbb C$$
in Hausdorff topology.
\end{thm}

\begin{pf}
Take any $z\notin\cup_{n<0} J({\cal T}_n)$. Then the orbit ${\cal O}_{\cal T}(z)$
escapes every domain $\Omega_n$ for large negative $n$.
That is, there exists an infinite sequence of points $z_{n_k}\in{\cal O}_{\cal T}(z)$
for $n_k\to -\infty$ such that $z_{n_k}\in\Omega_{n_k}$
and ${\cal H}_{n_k}(z_{n_k})\notin \Omega_{n_k}$.
Note, that by real {\it a priori} bounds and compactness of $\hol$, the difference
$|n_k-n_{k+1}|$ is bounded.
 By \lemref{expansion towers} (II) we have $||{\cal H}_{n_k}'(z_{n_k})||>C>1$ for some universal
constant $C$.
Moreover, as seen from the proof of the same Lemma, there is an element $h_{n_k}\in F({\cal T}^{{n_k}-1})$ 
 such that $\dist_{{\Bbb C}\setminus{\Bbb R}}(h_{n_k}(z_{n_k}),J({\cal T}^{{n_k}-1}))<s$
for some universal value of $s>0$. Let us arrange the points $\{\zeta_k=h_{n_k}(z_{n_k})\}$ into an 
infinite orbit
$\zeta_0=z,\zeta_1=g_0(\zeta_0),\zeta_2=g_1(\zeta_1)\ldots$,
 with $g_k\in F({\cal T})$. By 
\lemref{expansion towers} (I) we have $||g_k'(\zeta_k)||>C.$
Denote by $\alpha_k$ the hyperbolic geodesic in ${\Bbb C}\setminus{\Bbb R}$ of length $l(\alpha_k)<s$
connecting $\zeta_k$ to $J({\cal T}^{{n_k}-1})$.   Let $\alpha_{n_k}'$ be the connected component of
$g_0^{-1}\circ\ldots\circ g^{-1}_k(\alpha_{n_k})$ containing $z=\zeta_0$. 
Since the Julia set of a tower is invariant, it is enough to show that 
$l(\alpha'_{n_k})\to 0$ as $n_k\to-\infty$. Indeed, 
$$||D g_k\circ \cdots\circ g_0||>C^k.$$
By \lemref{variation distortion} this inequality holds along $\alpha'_{n_k}$ with $C$ replaced
by another universal constant $C_1>1$,
 and hence $\alpha_{n_k}'$ shrinks to $0$.
\end{pf}

Recall, that a {\it measurable line field} is a measurable Beltrami differential
$u(z)$ with $|u(z)|=1$ almost everywhere on the support of $u(z)$.
A line field corresponding to a Beltrami differential $\mu$ is 
$\mu/|\mu|$. We say that $u(z)$ is {\it invariant under a tower $\cal T$} if
for any $h\in F({\cal T})$, $h^*u/u$ is a real-valued function.
A measurable line field $u(z)$ is {\it univalent } on an open set $U$ if 
$u=h^*(d\bar z/dz)$  a.e. for a univalent map $h:U\to \Bbb C$.

\medskip
\noindent
{\it Proof of \thmref{bdl}}.
Suppose that ${\cal T}=({\cal H}_k)^\infty_{-\infty}$ is a bi-infinite limiting renormalization
tower, and $\mu$ is a nontrivial invariant Beltrami differential of $\cal T$. Let $u(z)$ denote
the corresponding invariant line field.
By \thmref{rigidity forward} $u(z)$ is not supported in $\cup J({\cal T}^n)$.
 Let $z$ be a point of almost continuity of $u(z)$, $z\notin J({\cal T}^n)$ for any $n$.
Denote again by $z_{n_k}\in \Omega_{n_k}$ the elements of ${\cal O}_{\cal T}(z)$ with
${\cal H}_{n_k}(z_{n_k})\notin \Omega_{n_k}$. 
As seen from the proof of  \lemref{expansion towers},
 there is an element $h_{n_k}\in F({\cal T}^{{n_k}-1})$ such that
the Euclidean distance from $\zeta_k=h_{n_k}(z_{n_k})$ to $\Bbb R$ is commensurable with
$|\Omega_{n_k}\cap {\Bbb R}|$.
Let ${\cal T}_{n_k}'=\Lambda_{n_k}\circ {\cal T}\circ \Lambda_{n_k}^{-1}$ where $\Lambda_{n_k}(z)=z/|\Omega_{n_k}\cap{\Bbb R}|$.
By compactness of $\hol$, we may ensure, passing to a subsequence, that elements of ${\cal T}_{n_k}'$
converge to holomorphic pairs ${\cal H}_i'$ forming a bi-infinite tower ${\cal T}'$.
By choosing a further subsequence we may assume that $\Lambda_{n_k}^{-1}\zeta_k\to w$
, and (see \cite{McM2})  rescaled linefields $u(\Lambda_{n_k}^{-1}(z))$ $w^*$-converge 
to a measurable linefield $u'$, with $h^*u'=u'$ for all $h\in F({\cal T'})$.

Let $D$ be a small disk around $w$ in ${\Bbb C}\setminus{\Bbb R}$, and denote by $D_{n_k}$ 
its image under the homothety $\Lambda_{n_k}$. The diameters of all $D_{n_k}$ in the hyperbolic
metric of ${\Bbb C}\setminus {\Bbb R}$ are equal. Denote by $D_{n_k}'$ the univalent preimage of
$D_{n_k}$ by maps of $F({\cal T})$ containing the point $z$.
\lemref{expansion towers} together with \lemref{variation distortion} imply 
that Euclidean diameters of $D_{n_k}$ shrink to $0$, and $z$ is well inside $D_{n_k}$. By 
\cite[Theorem 5.16]{McM1} we can choose the linefield $u'$ to be univalent on $D$.

By \thmref{fill the plane} there is an $i$ such that $J({\cal H}_i')\cap D\ne\emptyset. $
By \lemref{periodic dense} the orbit of $D$ by ${\cal T}'_i$ covers 
all of $\Delta'_i$. By invariance, this means that $u'$ coincides with a locally univalent line field
around $0$ and $({\cal H}_i')^2(0)$ which implies contradictory behaviour of $u'$ around ${\cal H}_i'(0)$.
$\;\;\Box$

\section{Conclusions}
\noindent
{\it Proof of Theorems A, B.} 
Due to compactness of $\hol$, for each bi-infinite string $\rho\in\bar\Sigma$ 
 there exits a bi-infinite limiting
renormalization tower ${\cal T}=({\cal H}_i)^\infty_{-\infty}$
with the rotation number $\rho$.
In view of \thmref{uniqueness} such a tower is unique. The element of the set $\cal A$ corresponding
to the string $\rho$ is set to be $\zeta_{\rho}={\cal H}_0\cap{\Bbb R}$.

Let $\zeta=\zeta_{(\ldots,r_{-n},\ldots,r_0,\ldots,r_n,\ldots)}
\in \cal A$ be a commuting pair with $r_0=\infty$. Let $m$ be the smallest positive integer
for which $r_m=\infty$ (if such a number does not exist, set $m=\infty$).
The generalized renormalization $\cal G$ of such a pair is defined as
$${\cal G}\zeta={\cal P}_\tau\zeta,$$
where the transit map $\tau$ is chosen so that 
$$\rho({\cal G}\zeta)=[r_1,\ldots,r_m],$$
and moreover, when $m$ is finite, ${\cal G}\zeta$ is parabolic.
By \lemref{one transit} such $\tau$ is unique. The action
of $\cal G$ on $\cal A$ is  injective (Remarks \ref{injective1},\ref{injective2});
the required invariance properties of $\cal A$ readily follow from \thmref{uniqueness}.

Let  $\zeta\in\cal E$ be any commuting pair with $\rho(\zeta)\in {\Bbb R}\setminus{\Bbb Q}$,
and let $\hat\zeta$ be a limit point of the sequence ${\cal R}^n\zeta$. By compactness of
$\hol$, $\hat\zeta$ is the base pair of a bi-infinite limiting renormalization tower
and therefore $\hat\zeta\in\cal A$. 

By compactness of $\hol$, there exists an open neighborhood $G$ of the origin, such that
the maximal domains of definition of the elements of ${\cal R}^n\zeta$ contain $G$ for all large $n$.
Now let $\zeta'\in \cal A$ be a commuting pair with $\rho(\zeta')=\rho(\zeta)$.
We will show that 
$$\dist({\cal R}^n\zeta,{\cal R}^n\zeta')\to 0$$
where the distance is measured as the maximum of the distance between the analytic extensions of the 
elements of the renormalized pairs in $C^0$-metric on $G$.
Indeed, otherwise there exists a sequence $n_k\to\infty$ and $\eps>0$ such that
$$\dist({\cal R}^{n_k}\zeta,{\cal R}^{n_k}\zeta')>\eps.$$
Passing to further subsequence we may ensure that ${\cal R}^{n_k}\zeta\to \zeta_1$,
${\cal R}^{n_k}\zeta'\to\zeta_2$ and $\zeta_1$ and $\zeta_2$ are the base pairs of two limiting renormalization
towers ${\cal T}_1$, ${\cal T}_2$ with $\rho({\cal T}_1)=\rho({\cal T}_2)$. 
The towers ${\cal T}_1$ and ${\cal T}_2$ are not homothetic, which contradicts \thmref{combinatorial rigidity}.
\hfill
$\Box$

\medskip
\noindent
{\it Proof of Theorem C.} For a parabolic pair $\zeta\in\cal E$ associate a string $(r_i)_1^\infty$,
$m\leq\infty$ with $r_i\in {\Bbb N}\cup\{\infty\}$. Inductively define
the pair ${\cal G}^{n+1}\zeta$ as ${\cal R}({\cal G}^n\zeta)$ when $\chi({\cal G}^n\zeta)\ne \infty$.
If ${\cal G}^n\zeta$ is a parabolic pair,
let $m$ be the smallest positive integer such that $r_{n+m}=\infty$ (if it does not exist, set $m=\infty$).
 Set ${\cal G}^{n+1}\zeta={\cal P}_\tau ({\cal G}^n\zeta)$ 
where the transit map $\tau$ is chosen so that $\rho({\cal G}^{n+1}\zeta)=[r_n,r_{n+1},\ldots,r_m]$
and, in the case when $m$ is finite, ${\cal G}^{n+1}\zeta$ is a parabolic pair (by \lemref{one transit}
such $\tau$ is unique). 

Consider a sequence of perturbations $\zeta_k\to \zeta$ with $\zeta_k\in\cal E$ and 
$$\rho(\zeta_k)=[r_0^k,r_1^k,\ldots]
\in{\Bbb R}\setminus{\Bbb Q}.$$ 
We will also require that $r_n^k\underset{k\to\infty}{\longrightarrow}r_n$.
In view of \lemref{ren-converge}, we may assume by choosing a further subsequence that  
$${\cal R}^n\zeta_k\underset{k\to\infty}{\longrightarrow}{\cal G}^n\zeta.$$
The existence of uniform real bounds for maps in the Epstein class implies that for $n>N_0$
the family
$\{{\cal R}^n\zeta_k\}$ is sequentially pre-compact in $\cal E$, which means that 
the  sequence $\{{\cal G}^n\zeta\}$ is pre-compact as well. 
Let $\zeta^*$ be any limit point of $\{{\cal G}^n\zeta\}$. For some choice of $n_i$ we have 
$${\cal R}^{n_i}\zeta_k\to\zeta^*.$$
Choosing a ``diagonal" subsequence $\zeta_{k_j}$ we may ensure that 
${\cal R}^{n_i+L}\zeta_{k_j}$ also converge for all $L\in \Bbb Z$. Thus $\zeta^*$ is a base map of 
a bi-infinite renormalization limiting tower, and hence $\zeta^*\in\cal A$. $\;\;\;\Box$

\medskip
\noindent
{\it Proof of Theorem D.} The pair $\zeta_0$ is the base pair of the bi-infinite renormalization
tower with rotation number $(\ldots,\infty,\infty,\infty,\ldots).$ The required properties of
$\zeta_0$ immediately follow from the above considerations.
$\;\;\;\Box$

\end{document}